\documentclass[11pt]{article}
\usepackage{amsmath,amsfonts,amsthm,epsf,stackrel}
\usepackage{graphicx,subfigure}
\usepackage{epstopdf}
\usepackage{multirow}
\usepackage{color}

\usepackage{slashbox}

\RequirePackage{natbib}

\usepackage[margin=2cm]{geometry}
\usepackage{mathabx}

\newcommand{\ignore}[1]{}
\newcommand{\tendp}{\stackrel{p}{\longrightarrow}}
\newcommand{\tendd}{\stackrel{{\cal D}}{\longrightarrow}}
\newcommand{\ve}{\varepsilon}
\newcommand{\bs}{\boldsymbol}

\newtheorem{theorem}{Theorem}
\newtheorem{corollary}{Corollary}
\newtheorem{proposition}{Proposition}
\newtheorem{lemma}{Lemma}

\begin{document}

\title{Semi-Supervised linear regression}
\author{David Azriel$^1$, Lawrence D. Brown$^2$, Michael Sklar$^3$, Richard Berk$^2$,\\
	 Andreas Buja$^2$ and Linda Zhao$^2$\thanks{Corresponding author: David Azriel, davidazr@technion.ac.il. Lawrence D. Brown deceased on February 21, 2018.}}
\date{}
\maketitle

\noindent{Technion -- Israel institute of Technology$^1$, Wharton -- University of Pennsylvania$^2$ and Stanford University$^3$}

\begin{abstract}
We study a regression problem where for some part of the data we observe both the {\em label} variable ($Y$) and the {\em predictors} (${\bf X}$), while for other part of the data only the predictors are given. Such a problem arises, for example, when observations of the label variable are costly and may require a skilled human agent. {When the conditional expectation $E[Y | {\bf X}]$ is not exactly linear, one can consider the best linear approximation to the conditional expectation, which can be estimated consistently by the least squares estimates (LSE). The latter depends only on the labeled data.   
We suggest improved alternative estimates to the LSE that use also the unlabeled data.} Our estimation method can be easily implemented and has simply described asymptotic properties.
The new estimates asymptotically dominate the usual standard procedures under certain non-linearity condition of $E[Y | {\bf X}]$; otherwise, they are asymptotically equivalent.
The performance of the new estimator for small sample size is investigated in an extensive simulation study. A real data example of inferring homeless population is used to illustrate the new methodology.     
\end{abstract}

\section{Introduction}

\subsection{Background and contribution}
The term ``semi-supervised learning'' was coined in the machine learning literature to describe a situation in which some of the data is labeled while the rest of the data is unlabeled \citep{Merz}. 
Such situations occur when the label variable is difficult to observe and may require a complicated or expensive procedure. A typical example is web document classification, where the classification is done by a human agent while there are many more unlabeled on-line documents. 
Specifically, 
a sample of $n$ observations from the joint distribution of $({\bf X},Y)$ is given, where $Y$ is a one-dimensional {\em label} variable and ${\bf X}$ is a p-dimensional vector of {\em covariates} or {\em predictors}. 
Also, an additional sample of size $m$ is observed where only ${\bf X}$ is given.
The purpose is to study procedures that make use of the additional unlabeled data to better capture the shape of the underlying joint distribution of the labeled data. 

A large body of literature focuses on the case that $Y$ takes a small number of values and the problem reduces to a classification task; see e.g., \citet{Zhu} and references therein. 
{When the predictors approximately lie in a low dimensional space, the unlabeled sample can be used to estimate the low dimensional structure  \citep[see e.g.,][]{Goldberg} }.
 \citet{Wang} divide the different methods into two approaches: distributional and margin-based.  
The distributional approach relies on an assumption relating the conditional expectation $E[Y | {\bf X}]$ to the marginal distribution of ${\bf X}$ and the margin-based approach uses the extra information on ${\bf X}$ for estimating the Bayes decision boundary; see \citet{Liang} for a Bayesian perspective within this basic approach.
Other works consider continuous $Y$'s and use the unlabeled data to learn the structure of the ${\bf X}$'s in order to better estimate a non-parametric regression \citep{Zhou,Lafferty,Johnson}. These works could be very helpful in situations where non-parametric regression is useful and unlabeled data are available. 





{
Here we follow a different methodology. We aim to estimate the vector ${\boldsymbol \beta}$ composed of the parameters of the best linear predictor of $Y$, but we do not necessarily assume that $E[Y | {\bf X}]$ is exactly linear in ${\bf X}$. The methodology uses possible unmodeled non-linearities in $E(Y|{\bf X})$ and adapts information from the unlabeled data to provide an estimator for ${\boldsymbol \beta}$. This estimator is asymptotically superior to the standard least squares estimate (LSE), both in terms of variance of ${\boldsymbol \beta}$ and mean squared error in predicting $Y$.}



{ In the statistical literature, the typical approach to regression with unlabeled data} may be best summarized by the following quote from \citet{Little}:
\begin{quote}
The related problem of missing values in the outcome $Y$ was prominent in the early history of missing-data methods, but is less interesting in the following sense: If the $X$'s are complete and the missing values of $Y$ are missing at random, then the incomplete cases contribute no information to the regression of $Y$ on $X_1, \ldots, X_p$.
\end{quote}
But see \citet[Chapter 7]{Cochran} for a different view more closely in tune with our development.

\citet{Buja} show that ${\boldsymbol \beta}$ does not depend on the distribution of ${\bf X}$ if and only if $E[Y | {\bf X}]$ is linear in ${\bf X}${, that is ${\bf X}$ is ancillary if and only if the linear model holds true. Therefore, by the conditionality principle \citep[see, e.g.,][Chapter 2.3]{Cox} the true marginal distribution of ${\bf X}$ is not useful for inferring ${\bs \beta}$. However, several counterexamples to this principle in the context of linear regression has been noted in the literature. 
\citet{Brown} constructs an estimate that depends on the marginal distribution of ${\bf X}$ and uniformly improve over the LSE even when ${\bf X}$ is conditioned upon. On the other hand, the improvement vanishes as the sample size goes to infinity. \citet{Ryan} study semi-supervised linear regression when the distribution of ${\bf X}$ in the unlabeled part is different from the labeled data. They suggested an estimator that accounts for this situation and show that it has a smaller mean square error     
 when a certain geometrical condition is satisfied. \citet{Azriel} argues that under a high-dimensional linear regression model, consistent estimation of $Var(Y|{\bf X})$ is possible only when ${\bf X}$ is not conditioned upon, which raises questions about the applicability of the conditionality principle in this context.} 
 
The assumption that $E[Y | {\bf X}]$ is linear is unrealistic in many situations, and we show that in the absence of such an assumption the unlabeled sample can be used to provide useful information for the estimation of ${\boldsymbol \beta}$ and uniform improvement over the LSE can be achieved (asymptotically).   

We consider two scenarios. {In the first, the distribution of ${\bf X}$ is known exactly. This is equivalent to having infinitely many unlabeled observations ($m=\infty$). We call this the \textit{full information} scenario. In the second scenario, which is more frequently encountered, $m < \infty$ but we assume that $m$ is at least proportional to n.  We call this \textit{partial information}. In both situations we provide an asymptotically better estimate of ${\boldsymbol \beta}$ than the standard LSE, and hence also a better linear predictor. Our new procedure is closely related to the semi-supervised estimations of means, which is the topic of a recent manuscript \citep{Anru}. This latter work constructs a semi-supervised regression estimator that is better than the sample mean. We build upon this work in two ways. First, we transform the regression problem into a mean estimation problem for each parameter via a ``multiplication step''. 
 Second, although it is proven by Zhang et al. that the semi-supervised estimate is superior to the sample mean, this does not directly imply superiority of our estimator to the LSE. For this purpose, a new argument is introduced based on an orthogonal decomposition of the error (see Section \ref{sec:nutshell} below).}

{In a recent paper, \citet{Chakrabortty} also construct an estimator that is asymptotically superior to the LSE. Chakrabortty and Cai estimate the non-linear part of $E(Y|{\bf X})$, non-parametrically, and then impute the missing $Y$'s to aid estimation of ${\boldsymbol \beta}$. This is challenging, and they suggest several methods for the nonparametric estimation.} 

{Our approach, by contrast, squeezes extra efficiency from the LSE using minimal modeling. The two approaches are discussed further in Section \ref{sec:simple}, where we show that our estimator is the ``simplest estimator'' that dominates the LSE in some sense. In a simulation study (Section \ref{sec:comp_sim}), we find that the suggested method is especially useful when the non-linear part is difficult to estimate for the sample size at hand.}
	
%
%

\subsection{Assumption lean linear regression}\label{sec:assumption}

A motivating premise of this research is that linear regression can and should be used even in the presence of non-linearities. Recalling G.E.P. Box’s famous dictum, we believe the linear model is often useful even when linearity is unverifiable or wrong. \citet{Buja} and \citet{Berk1}  discuss the use of misspecified linear regression at length; here we briefly mention some key conclusions:

First, the parameters of the linear model can be interpreted as functionals of the joint distribution of $({\bf X},Y)$. The typical interpretation of how a unit change in ${\bf X}$ affects the conditional distribution of $Y$ given ${\bf X}$ is not justified, but remains true in an averaged sense. The model fit may be interpreted as the best linear approximation to the data.

Second, the definition of the model parameters depends on the distribution of the vector of covariates ${\bf X}$, and therefore inference should not be conditioned on the observed covariates. In fact, in this current work, we find that information from the unlabeled data about the marginal distribution of the covariates is useful for estimation.

Third, the standard variance of the LSE is wrong and the ``sandwich'' form (see \eqref{eq:sandV} below) should be used instead. The variance could be consistently estimated either by using a parametric estimator or by the bootstrap; see Section \ref{sec:est_se} and Appendix \ref{sec:app_est}.   

There are several reasons to use the linear model even when it is wrong. For a given dataset, a data-driven model selection procedure or trying out different transformations, may lead to over-fitting. Furthermore, some covariates may express phenomena in ``natural'' or ``conventional'' units that should not be transformed even if model fit is improved. Also, the linear model provides a meaningful measure of the association between the response and the covariates that is useful in many cases. Finally, for small sample-sizes, the linear model can provide in certain datasets a smaller prediction error than non-parametric alternatives. For a more elaborate discussion of the usefulness of the linear model see \citet{Buja} and \citet{Berk1}.       

\subsection{{The suggested methodology in a nutshell}}\label{sec:nutshell}
To demonstrate our basic methodology, consider the simple linear regression model,
\begin{equation} \label{eq:model_1}
Y = \alpha + \beta X + \delta,
\end{equation}
where $\alpha$, $\beta$ {are least squares coefficients and $\delta$ is a remainder term; these will be defined in more detail shortly}. When $E(X)$ and $E(X^2)$ are known, we construct an estimator of $\beta$ that asymptotically dominates the least squares estimator. This is done by replacing model \eqref{eq:model_1} with a different model where $\beta$ is the intercept {and also the expectation of the newly defined label variable.  For this latter model, \citet{Anru} show that the intercept estimator dominates the simple empirical estimate in a semi-supervised setting. {Here we go one step further} and show that the intercept estimator in this new model is better than the least squares estimator from the original model as explained below.}

To present the intercept model and the methodology,
begin by stating the general form of the best linear approximation. The approximation is used in two different models, and therefore we introduce now the notation ${\bf U}$ and $W$ instead of ${\bf X}$ and $Y$.
Suppose that $W \in \mathbb{R}$ and ${\bf U} \in \mathbb{R}^d$ are random variables with joint distribution $G$ and finite second moments. The best linear predictor is 
\begin{equation} \label{eq:best_lin}
{\boldsymbol \theta}=({\theta}_1,\ldots,{\theta}_d)^T={\arg \min}_{\tilde {\boldsymbol \theta} \in {\mathbb R}^d } E ( W -\sum_{j=1}^d U_j \tilde{\theta}_j )^2 = \left\{ E( {\bf U} {\bf U}^T ) \right\}^{-1} E( {\bf U} W ). 
\end{equation}    
Notice that this is a population version of the least squares {where $\sum_{j=1}^d U_j \tilde{\theta}_j$ minimizes the $L_2$ distance in the population from $W$ to any linear function of ${\bf U}$.}
It follows that
\begin{equation} \label{eq:best_reg}
W = \sum_{j=1}^d \theta_j U_j + r,
\end{equation}
where the remainder term $r = W - \sum_{j=1}^d \theta_j U_j$ is orthogonal to ${\bf U}$, i.e., $E(r {\bf U})=0$. 
Given a sample of $n$ observations from  $G$, the standard least squares estimate, $\hat{\boldsymbol \theta}_{LSE}$, 
satisfies asymptotically  \citep{White1}
\begin{equation} \label{eq:sandV}
\sqrt{n}( \hat{\boldsymbol \theta}_{LSE} - {\boldsymbol \theta}) \tendd N\left( 0 , \left\{ E({\bf U} {\bf U}^T ) \right\}^{-1} E\left(r^2  {\bf U} {\bf U}^T  \right) \left\{ E( {\bf U} {\bf U}^T ) \right\}^{-1} \right).
\end{equation}
{Unlike the standard fixed ${\bf X}$ assumption, the asymptotic variance has a ``sandwich'' form, $\left\{ E({\bf U} {\bf U}^T ) \right\}^{-1}$ forming the ``bread'' and $E\left(r^2  {\bf U} {\bf U}^T  \right)$ the ``meat''.} {See \citet{Buja} for further discussion of this form of the sandwich.}

We now return to model \eqref{eq:model_1} and assume that $E(X)$ and $E(X^2)$ are known. We can therefore assume w.l.o.g that $X$ is standardized (i.e., $E(X)=0$ and $E(X^2)=1$). In this case we can write model \eqref{eq:model_1} with $\alpha = E(Y)$, $\beta=E(XY)$ and $\delta = Y-(\alpha + \beta X)$; here we consider model \eqref{eq:best_reg} with $W=Y$, $U_1=1$, $U_2=X$ and the remainder term is $r=\delta$. 
The standard LSE is a consistent estimate for $\beta$ and satisfies, according to \eqref{eq:sandV},
\[
\sqrt{n} (  \hat{\beta}_{LSE} - \beta) \tendd N\left\{ 0, E( \delta X)^2\right\}. 
\]  

{Our aim it to better estimate $\beta=E(XY)$. To this end, we multiply \eqref{eq:model_1} by $X$ setting $XY$ to be the labeled variable $W$. Furthermore, we also center the varibale ${\bf U}$ setting $E(XY)=\beta$ to be the intercept. Specifically, }
\begin{equation} \label{eq:model_2}
XY = \alpha X + \beta X^2  +  X \delta = \beta + a X + b (X^2 -1) + \tilde{\delta}.
\end{equation}  
Now we consider model \eqref{eq:best_reg} with $W=XY$, $U_1=1$, $U_2=X$, $U_3=X^2-1$ and $r=\tilde{\delta}$; here $a,b$ are $\theta_2,\theta_3$ defined by \eqref{eq:best_lin}. In setting $\beta=E(XY)$ to be the intercept coefficient we used that $EX=E(X^2)-1=0$, i.e., that we know the first and second moment of $X$.
We define our new estimator, $\hat{\beta}_{TI}$ (TI for total information), to be the {LSE intercept estimator of \eqref{eq:model_2} based on $n$ observations}. 
The sandwich theorem \eqref{eq:sandV} implies that 
\[
\sqrt{n} (  \hat{\beta}_{TI} - \beta) \tendd N( 0, E \tilde{ \delta} ^2 ). 
\]  

Since $\delta X - \tilde{\delta}= X(a-\alpha)+(X^2-1)(b-\beta)$ is a linear function of $X,X^2$, then {it is orthogonal to the remainder term of \eqref{eq:model_2}, which is $\tilde{\delta}$, (i.e., $E(\delta X - \tilde{\delta})\tilde{\delta}=0$)} and therefore,
\begin{equation}\label{eq:decomp}
E( \delta X)^2 = E\{ \tilde{\delta} + (\delta X - \tilde{\delta})\}^2= E\tilde{\delta}^2 + E(\delta X - \tilde{\delta})^2.  
\end{equation}
This implies that the asymptotic variance of $\hat{\beta}_{TI}$ is smaller than that of $\hat{\beta}_{LSE}$ with equality iff $\delta X \equiv \tilde{\delta}$. The latter occurs iff $\alpha=a$ and $b=\beta$ or equivalently that $\delta$ is uncorrelated with both $X^2$ and $X^3$. {This occurs when the non-linear part of $E[Y|X]$ (if exists) is uncorrelated with $X^2$ and $X^3$.}
 {In this case, Models \eqref{eq:model_1} and 
\eqref{eq:model_2} are essentially the same model and nothing is gained in the new methodology. On the other hand, when $\delta$ is correlated with either $X^2$ or $X^3$, then Models \eqref{eq:model_1} and \eqref{eq:model_2} are different and $\hat{\beta}_{TI}$ has smaller asymptotic variance than $\hat{\beta}_{LSE}$ (see Appendix \ref{app:F}).
We further show that a similar decomposition to \eqref{eq:decomp} holds for the partial information case and we generalize the results to the $p$-dimensional case. {We also show that smaller variance implies smaller prediction error and therefore $\hat{\beta}_{TI}$ provides better predictions than $\hat{\beta}_{LSE}$.}

The rest of the paper is organized as follows. 
Section \ref{sec:prel} provides the basic setting and the loss functions that we use. The main results are given in Section \ref{sec:main} and in Section \ref{sec:simple} we discuss the relation between our method and semi-parametric efficiency.  Section \ref{sec:sim} describes an extensive simulation study and 
Section \ref{sec:est_se} discusses estimation of the asymptotic covariance matrix of the estimates. An
implementation of the new methodology to infer homeless population in Los Angeles is discussed in Section \ref{sec:inferring}. Section \ref{sec:dis} concludes with final remarks. The proofs are given in Section \ref{sec:proofs}.

\section{Preliminaries} \label{sec:prel}

Consider a sample of $n$ i.i.d observations $({\bf X}^{(1)},Y^{(1)}),\ldots,({\bf X}^{(n)},Y^{(n)})$ from a joint distribution $G$, where ${\bf X} \in \mathbb{R}^p$, $Y \in \mathbb{R}$, and an additional set of $m$ independent observations $({\bf X}^{(n+1)}\ldots,{\bf X}^{(n+m)})$ from the marginal distribution of ${\bf X}$. We use super-index to denote the number of the observation, and sub-index to denote coordinates of ${\bf X}$. The notation ${\bf X}, Y$ without super-index denotes a random vector whose distribution is $G$.

We write $\vec{\bf X}=(1,X_1,\ldots,X_p)^T$ to be a vector ${\bf X}$ with an additional constant 1 to accommodate an intercept term.  
Assume that the second moments of $G$ exists and that the matrix $E\left( \vec{\bf X} \vec{\bf X}^T \right)$ is invertible.           
Then, we can define 
\begin{equation}\label{eq:def_beta}
(\alpha, {\boldsymbol \beta}) = {\arg \min}_{\tilde{\alpha} \in \mathbb{R}, \tilde{\boldsymbol \beta} \in \mathbb{R}^p} E(Y - \tilde{\alpha} - \tilde{\boldsymbol \beta}^T X)^2 =  \left\{ E\left( \vec{\bf X} \vec{\bf X}^T \right) \right\}^{-1} E\left( \vec{\bf X} Y\right). 
\end{equation}

In the presence of non-linearity, ${\boldsymbol \beta}$ is still a meaningful parameter that describes the overall association between $Y$ and ${\bf X}$ \citep{Buja}.  
We have in mind two related purposes. The first purpose is just to better estimate the parameters of interest, while the second purpose pertains to prediction. The latter is formalized now.
Suppose that an independent observation $({\bf X}^*,Y^*) \sim G$ is given. The optimal linear predictor is $\alpha + {\boldsymbol \beta}^T X^*$. 
We consider the excess loss of an estimator $\tilde{\alpha},\tilde{\boldsymbol \beta}$ 
\[
L(\tilde{\alpha},\tilde{\boldsymbol \beta})= \left(Y^* - \tilde{\alpha} - \tilde{\boldsymbol \beta}^T {\bf X}^*\right)^2- \left(Y^* - {\alpha} - {\boldsymbol \beta}^T {\bf X}^* \right)^2.
\]
We have that (by Lemma \ref{lem1})
\begin{equation*}
 E L(\tilde{\alpha},\tilde{\boldsymbol \beta})=  E\Big[  \bar{Y} - \sum_{j=1}^p \tilde{\beta}_j\{ {\bar X}_j - E(X_j)\} -E(Y)\Big]^2 +  E\left\{ ({\boldsymbol \beta}- \tilde{\boldsymbol \beta})^T {\bf M} ({\boldsymbol \beta}- \tilde{\boldsymbol \beta}) \right \}= E \tilde{L}(\tilde{\alpha},\tilde{\boldsymbol \beta}),
\end{equation*}
where $\tilde{L}$ is the expression inside the expectation and ${\bf M}$ is the covariance matrix of ${\bf X}$. {Notice that $L$ and $\tilde{L}$ have the same expectation but they are different random variables.}
We further define the asymptotic risk as
\[
R(\tilde{\alpha},\tilde{\boldsymbol \beta}) = \lim_{B \to \infty} \lim_{n \to \infty}  E \min\{ n\tilde{L}(\tilde{\alpha},\tilde{\boldsymbol \beta}), B\}.
\]
{The loss is of order $1/n$ and therefore we consider expectation of $n\tilde{L}$. This is truncated by {an arbitrarily large number $B$}, since when the loss is large, it makes sense not to penalize any further.} Also, the truncation helps to avoid issues of uniform integrability. This is done for example in \citet{Le_Cam}, Chapter 5.
The following proposition provides a simple expression for the asymptotic risk.
\begin{proposition}\label{prop1}
Let  $\tilde{\alpha},\tilde{\boldsymbol \beta}$ satisfies 
\begin{equation}\label{eq:tilde_alph}
\tilde{\alpha} = \bar{Y} - \sum_{j=1}^p \tilde{\beta}_j {\bar X}_j,
\end{equation}
and assume that $\tilde{\boldsymbol \beta}$ satisfies $\sqrt{n} (\tilde{\boldsymbol \beta}- {\boldsymbol \beta}) \tendd N(0 , {\boldsymbol \Sigma})$ then,
\begin{equation}\label{eq:excess_risk1}
{R}(\tilde{\alpha},\tilde{\boldsymbol \beta}) =  E\Big[  Y - \sum_{j=1}^p {\beta}_j\{ { X}_j - E(X_j)\} -E(Y)\Big]^2
+ {\rm Trace} ({\bf M} {\boldsymbol \Sigma}),
\end{equation}
where ${\bf M}$ is the covariance matrix of ${\bf X}$.
\end{proposition}
The first term in \eqref{eq:excess_risk1} does not depend on the distribution of $\tilde{\boldsymbol \beta}$. Hence,
Proposition \ref{prop1} shows that the excess risk is minimized when $\text{Trace} ({\bf M} {\boldsymbol \Sigma})$ is small. Thus, we aim at estimators $\tilde{\boldsymbol \beta}$  with asymptotic distribution $N(0,{\boldsymbol \Sigma})$ such that ${\boldsymbol \Sigma}$ is ``smaller'' than the covariance matrix of LSE, in the sense that the difference is positive semidefinite. Such an estimator asymptotically better estimates ${\boldsymbol \beta}$ and also has smaller asymptotic excess risk. 

For two estimates $(\tilde{\alpha}^{(1)},\tilde{\bs \beta}^{(1)})$ and $(\tilde{\alpha}^{(2)},\tilde{\bs \beta}^{(2)})$, Proposition \ref{prop1} implies that
\[
E L(\tilde{\alpha}^{(1)},\tilde{\boldsymbol \beta}^{(1)}) - E L(\tilde{\alpha}^{(2)},\tilde{\boldsymbol \beta}^{(2)}) \approx \text{Trace} \left\{ {\bf M}( {\boldsymbol \Sigma}^{(1)}-{\boldsymbol \Sigma}^{(2)} )\right\} / n,
\]
and therefore the difference of the prediction errors is
\begin{equation}\label{eq:pred_error1}
E\left\{Y^* - \tilde{\alpha}^{(1)} - \left(\tilde{\boldsymbol \beta}^{(1)} \right)^T {\bf X}^*\right\}^2- E\left\{Y^* - \tilde{\alpha}^{(2)} - \left(\tilde{\boldsymbol \beta}^{(2)} \right)^T {\bf X}^*\right\}^2 \approx \text{Trace} \left\{ {\bf M}( {\boldsymbol \Sigma}^{(1)}-{\boldsymbol \Sigma}^{(2)} )\right\} / n.
\end{equation}
It follows that the difference of the errors of the prediction of the mean is also
\begin{equation}\label{eq:pred_error2}
\small{
E\left\{ E(Y^*|{\bf X}^*) - \tilde{\alpha}^{(1)} - \left(\tilde{\boldsymbol \beta}^{(1)} \right)^T {\bf X}^*\right\}^2- E\left\{E(Y^*|{\bf X}^*) - \tilde{\alpha}^{(2)} - \left(\tilde{\boldsymbol \beta}^{(2)} \right)^T {\bf X}^*\right\}^2 \approx \text{Trace} \left\{ {\bf M}( {\boldsymbol \Sigma}^{(1)}-{\boldsymbol \Sigma}^{(2)} )\right\}/n}.
\end{equation}

Notice that we estimate the intercept, as in \eqref{eq:tilde_alph}, using ${\bar X}_j$ rather than the possibly known expectations $E(X_j)$. This is in accordance with the findings of \citet{Anru}. 

\section{Main results} \label{sec:main}
 
In this section we provide the main theoretical results of the paper. 
Before we deal with the p-dimensional case, we first introduce the results for one-dimensional ${\bf X}$. The reason is twofold: first, the presentation in the one-dimensional case is simpler and captures the main ideas, and second, our results for the p-dimensional ${\bf X}$ are obtained by reducing the regression problem to p problems, each of which is closely related to the one dimensional ${\bf X}$ regression.  

\subsection{One dimensional ${X}$}\label{sec:one_dim}

{In this section we study the one dimensional case. Summary of the notation used here is presented in Table \ref{tab:not1}.} When ${G}$ has finite second moments we can write
\begin{equation}\label{eq:model1_1}
Y= \alpha + \beta X + \delta,
\end{equation}
where $\alpha = E(Y) - \beta E(X)$, $\beta = \frac{{\rm Cov}(X,Y)}{Var(X)}$ and $\delta=Y-\alpha - \beta X$. {Equation \eqref{eq:model1_1} is the best linear approximation in the population in the sense that} $\alpha + \beta X$ minimizes $E(\tilde{\alpha} + \tilde{\beta} X - Y)^2$ over all $\tilde{\alpha},\tilde{\beta} \in \mathbb{R}$ and $\delta$ is orthogonal to $1,X$, i.e., $E(\delta)=E(\delta X) = 0$.
The regular LSE is 
\[
\hat{\beta}_{LSE} = \frac{\sum_{i=1}^n (X^{(i)}-\bar{X}) Y^{(i)}}{\sum_{i=1}^n (X^{(i)}-\bar{X})^2}.
\]
This is a consistent and asymptotically unbiased and asymptotically normal estimator for $\beta$. See \citet{Buja} for contemporary discussion.

\begin{table} [ht!] 
\caption{\footnotesize Summary of the notation used in Section \ref{sec:one_dim}. }
\begin{center}
{
  \begin{tabular}{c c c}
    \hline
		Basic model &  $Y= \alpha + \beta X + \delta$ &  $\alpha = E(Y) - \beta E(X)$, $\beta = \frac{{\rm Cov}(X,Y)}{Var(X)}$ and $\delta=Y-\alpha - \beta X$\\
		\hline
    Intercept model & $W = \beta + a U_1  + b U_2  + \tilde{\delta}$ & $W=\frac{Y\{X-E(X)\}}{Var(X)}$, $U_1=\frac{X-E(X)}{Var(X)}$, $U_2=\frac{\{X-E(X)\}X}{Var(X)} -1$\\
		{} & {} & $a,b$ are $\theta_2,\theta_3$  defined by \eqref{eq:best_lin} for model \eqref{eq:modelXY}  \\
		{} & {} & and $\tilde{\delta}$ is the remainder term\\
		\hline
		TI estimator &  $\hat{\beta}_{TI} =  \bar{W}   - \hat{a}\bar{U}_1 -\hat b \bar{U}_2$ & $\bar{W} = \frac{1}{n} \sum_{i=1}^n W^{(i)}$, $\bar{U}_1=\frac{1}{n} \sum_{i=1}^n U_1^{(i)}$, $\bar{U}_2=\frac{1}{n} \sum_{i=1}^n U_2^{(i)}$\\ 
		{} & {} & $\hat{a},\hat{b}$ are the LSE of model \eqref{eq:modelXY}\\
		\hline
		PI estimator & $\hat{\beta}_{PI} = \bar{\check{W}}   - \hat{a}\bar{\check{U}}_1 -\hat b \bar{\check{U}}_2$ & $\bar{\check{W}}, \bar{\check{U}}_1, \bar{\check{U}}_2$ are the means over the labeled sample,\\
		{} & {} & 
		$\check{W}^{(i)}= \frac{\{X^{(i)}-\check{E}(X)\}Y^{(i)}}{\widecheck{Var}(X)}, \check{U}_1^{(i)}=\frac{\{X^{(i)} - \check{E}(X)\}}{\widecheck{Var}(X)}$\\
		{} & {} & $\check{U}_2^{(i)}=\frac{\{X^{(i)} - \check{E}(X)\}X^{(i)}}{\widecheck{Var}(X)}-1$, $\hat{a},\hat{b}$ are the LSE\\
	 {} & {} &  of model \eqref{eq:modelXY} with $\check{W},\check{U}_1,\check{U}_2$ instead of $W,U_1,U_2$\\
	\hline
	Asymptotics & $\nu =\lim \frac{n}{n+m} $ & $n$ ($m$) the labeled (unlabeled) sample size \\
	{} & $\sigma^2_{LSE}$, $\sigma^2_{PI}$, $\sigma^2_{TI}$ & asymptotic variance of the LSE, PI, TI estimators \\
	{} & {} & as in Theorem \ref{thm:one_dim}\\
	\hline 
	Main results {}&  $\sigma^2_{LSE} =  \sigma^2_{TI}+\sigma^2_{\rm diff}$ & $\sigma^2_{\rm diff}=E\left[ \frac{ \delta \{X - E(X)\}}{Var(X)}-\tilde{\delta}\right]^2$\\
	{} & $\sigma^2_{PI}= \sigma^2_{TI}+\nu\sigma^2_{\rm diff}$ & {} \\
	\hline
 		\end{tabular}}
\label{tab:not1}
\end{center}
\end{table}

For the total information estimator (TI) we consider the following regression model, which is obtained by multiplying \eqref{eq:model1_1} by $\frac{X-E(X)}{Var(X)}$,
\begin{equation} \label{eq:modelXY}
W = \beta + a U_1  + b U_2  + \tilde{\delta},  
\end{equation}
where $W=\frac{Y\{X-E(X)\}}{Var(X)}$, $U_1=\frac{X-E(X)}{Var(X)}$, $U_2=\frac{\{X-E(X)\}X}{Var(X)} -1$,
 $a,b$ are the coefficients $\theta_2,\theta_3$ of the best linear approximation defined by \eqref{eq:best_lin} and $\tilde{\delta}$ is the remainder term. {Here again $\beta + a U_1  + b U_2$ is the best linear approximation of $W$ in the population. The multiplication term makes the expectation of $W$ to be $E\left[\frac{Y\{X-E(X)\}}{Var(X)}\right]=\frac{Cov(X,Y)}{Var(X)}=\beta$. } 
The covariates $U_1, U_2$ are the predictors $1,X$ in \eqref{eq:model1_1} multiplied by  $\frac{X-E(X)}{Var(X)}$. 
We subtract 1 in $U_2$ in order to set $E(U_2)=0$. Since $E(U_1)=E(U_2)=0$, then the intercept is $E(W)= \beta$. Thus,   
we define the total information estimator to be the intercept estimator of \eqref{eq:modelXY}, i.e., 
\begin{equation}\label{eq:TI_CF}
\hat{\beta}_{TI} =  \bar{W}   - \hat{a}\bar{U}_1 -\hat b \bar{U}_2,
\end{equation}
where $\hat{a},\hat{b}$ are the regular LSE of model \eqref{eq:modelXY}, and $\bar{\cdot}$ denotes the mean over the supervised sample with $n$ observations. 

For the partial information (PI) estimator we estimate $W, U_1, U_2$ as follows: 
\[
\check{W}^{(i)}= \frac{\{X^{(i)}-\check{E}(X)\}Y^{(i)}}{\widecheck{Var}(X)}, \check{U}_1^{(i)}=\frac{\{X^{(i)} - \check{E}(X)\}}{\widecheck{Var}(X)}\text{ and }\check{U}_2^{(i)}=\frac{\{X^{(i)} - \check{E}(X)\}X^{(i)}}{\widecheck{Var}(X)}-1~ ,i=1,\ldots,n,
\]
where
$\check{E}(X) = \frac{1}{n+m}\sum_{i=1}^{n+m} X^{(i)}$,  $\widecheck{Var}(X)=\check{E}(X^2)-\{\check{E}(X)\}^2$, and
 $\check{E}(X^2) = \frac{1}{n+m}\sum_{i=1}^{n+m} \{X^{(i)}\}^2$.
The partial information estimator is
\[
\hat{\beta}_{PI} = \bar{\check{W}}   - \hat{a}\bar{\check{U}}_1 -\hat b \bar{\check{U}}_2,
\]
where $\hat{a}, \hat{b}$ are the regular LSE of the regression model \eqref{eq:modelXY} with $\check{W},\check{U}_1,\check{U}_2$ instead of $W,U_1,U_2$. (The estimates $\hat{a}, \hat{b}$ are different in the total and partial information cases but the same notation is presented for simplicity.) 
We use $\bar{\cdot}$ (respectively, $\check{\cdot}$) to denote empirical mean with respect to the labeled $n$ (respectively, full $n+m$) sample.
The following theorem states the asymptotic distribution of $\hat{\beta}_{LSE},\hat{\beta}_{TI}$ and $\hat{\beta}_{PI}$. The first part of the theorem is known and is stated here for comparison purposes. 

\begin{theorem} \label{thm:one_dim}
{}
\begin{enumerate}
\item[(i)] Suppose that $Var(X)\in (0 , \infty)$ and that $\delta \{X-E(X)\}$ has finite second moment, then,
\[
\sqrt{n}(\hat{\beta}_{LSE}-\beta) \tendd N(0,\sigma^2_{LSE}),
\]
where 
$\sigma^2_{LSE}= E\left[ \frac{ \delta^2 \{X - E(X)\}^2}{\{Var(X)\}^2}\right]$.
\item[(ii)]  Suppose further that the vector $(W,U_1,U_2)$ has finite second moments and that the matrix $E({\bf U} {\bf U}^T)$ is invertible for ${\bf U}= \left( 1, U_1, U_2\right)^T$.
Then,
\[
\sqrt{n}(\hat{\beta}_{TI}-\beta) \tendd N(0,\sigma^2_{TI}),
\]
where $\sigma^2_{TI}= E(\tilde{\delta}^2)$ and $\sigma^2_{LSE} =  \sigma^2_{TI}+\sigma^2_{\rm diff}$ where 
$\sigma^2_{\rm diff}=E\left[ \frac{ \delta \{X - E(X)\}}{Var(X)}-\tilde{\delta}\right]^2$.
\item[(iii)] Suppose further that $\lim \frac{n}{n+m} = \nu$, then, 
\[
\sqrt{n}(\hat{\beta}_{PI}-\beta) \tendd N(0,\sigma^2_{PI}),
\]
where $\sigma^2_{PI}= \sigma^2_{TI}+\nu\sigma^2_{\rm diff}$.
\end{enumerate}
Therefore, if $\sigma^2_{\rm diff} >0$ then $\sigma^2_{TI}<\sigma^2_{LSE}$ and if further $\nu<1$ then $\sigma^2_{PI}<\sigma^2_{LSE}$.
\end{theorem}
\begin{corollary} \label{cor:one_dim}
Theorem \ref{thm:one_dim} and Proposition \ref{prop1} imply that $R(\hat{\alpha}_{LSE},\hat{\beta}_{LSE})- R(\hat{\alpha}_{PI},\hat{\beta}_{PI}) = (1-\nu)\sigma^2_{\rm diff} Var(X)$.
\end{corollary}

The improvement of $\hat{\beta}_{TI}$ and $\hat{\beta}_{PI}$ over $\hat{\beta}_{LSE}$ depends on the assumption that $\sigma^2_{\rm diff}>0$. The quantity $\sigma^2_{\rm diff}$ measures the difference between the original regression model \eqref{eq:model1_1} and the intercept model \eqref{eq:modelXY}. If $E(\delta X^2)=E(\delta X^3)=0$ then $a=\alpha$ and $b=\beta$ and the two models are essentially the same, in which case $\sigma^2_{\rm diff}=0$, otherwise $\sigma^2_{\rm diff}>0$ (see Appendix \ref{app:F}). {In other words, when $E[Y|X]$ is non-linear and the non-linear part is correlated with either $X^2$ or $X^3$ then $\sigma^2_{\rm diff}>0$.}

 \subsection{Multidimensional ${\bf X}$} \label{sec:mult}
  
	We now consider the general $p$-dimensional case as described in the beginning of Section \ref{sec:prel}. {The notation of this section is summarized in Table \ref{tab:not2}.} The model can be written as
\begin{equation}\label{eq:model}
Y= \alpha + \beta_1 X_1+ \cdots +\beta_p X_p + \delta,
\end{equation}
where $\alpha, {\boldsymbol \beta}$ are the coefficients of the best linear predictor in the population defined by \eqref{eq:def_beta}. The remainder term $\delta=Y-\sum_{j=1}^p \beta_j X_j$ satisfies $E(\delta)=E(\delta X_1)=\cdots=E(\delta X_p)=0$. Our aim here is at estimating ${\boldsymbol \beta}$.
We use the adjustment representation of \citet{Buja} that reduces the p-dimensional estimation procedure to $p$ separate simple regression problems. Correspondingly, we will define our new estimates by solving $p$ mean-estimation-problems separately, one for each coordinate $j$.

\begin{table} [ht!] 
\caption{\footnotesize Summary of the notation used in Section \ref{sec:mult} when $j$ is the chosen coordinate.}
\begin{center}
\small{
  \begin{tabular}{c c c}
    \hline
		Basic model &  $Y= \alpha + \sum_{j=1}^p \beta_j X_j + \delta$ &  $\alpha, {\boldsymbol \beta}$ are  defined by \eqref{eq:def_beta}\\
		\hline
		Adjusted regressor &  $X_{j\bullet} = X_j -  \vec{\bf X}_{-j}^T {\boldsymbol \beta}_{-j \bullet }$ &  ${\boldsymbol \beta}_{-j \bullet }= \left\{E\left( \vec{\bf X}_{-j} \vec{\bf X}_{-j}^T \right) \right\}^{-1} E\left( \vec{\bf X}_{-j} X_j \right)$\\
		{} & {} & $\vec{\bf X}_{-j} = (1,X_1,\ldots,X_{j-1},X_{j+1},\ldots,X_p)^T$\\
		\hline
    Intercept model & $W_j = \beta_j +a U_1 + \sum_{j'=1}^pb_{j'} U_{j'+1} + \tilde{\delta}_j$&  $W_j = \frac{YX_{j\bullet}}{E(X_{j\bullet}^2)},U_1 =\frac{X_{j\bullet}}{E(X_{j\bullet}^2)}$ \\
		{} & {} & $U_{j'+1}=\frac{X_{j'}X_{j\bullet}}{E(X_{j\bullet}^2)} \text{ for }j'\ne j$, $U_{j+1}=\frac{X_{j}X_{j\bullet}}{E(X_{j\bullet}^2)}-1 $,\\
		{} & {} & $a,{\bf b}$ are $\theta_2,\cdots,\theta_{p+1}$  defined\\ 
		{} & {} & by \eqref{eq:best_lin} for model \eqref{eq:model1},\\
		{} & {} & $\tilde{\delta}_j$ is the remainder term\\
		\hline
		TI estimator &  $\{\hat{\boldsymbol \beta}_{TI}\}_j = \bar{W}_j -\hat{a} \bar{U}_1 -\sum_{j'=1}^p\hat{b}_{j'} \bar{U}_{j'+1}$  & $\bar{\cdot}$ is the mean over the labeled sample\\ 
		{} & {} & $\hat{a},\hat{\bf b}$ are the LSE of model \eqref{eq:model1}\\
		\hline
		PI estimator & $\{\hat{\boldsymbol \beta}_{PI}\}_j = \bar{\check{W}}_j -\hat{a} \bar{\check{U}}_1 -\sum_{j'=1}^p \hat{b}_{j'} \bar{\check{U}}_{j'+1}$ & $\bar{\cdot}$ is the mean over the labeled sample\\
		{} & {} & 
		$\check{W}_j = \frac{YX_{j\check{\bullet}}}{\check{E}(X_{j\check{\bullet}}^2)},~\check{U}_1 =\frac{X_{j\check{\bullet}}}{\check{E}(X_{j\check{\bullet}}^2)}$\\
		{} & {} & $\check{U}_{j'+1}=\frac{X_{j'}X_{j\check{\bullet}}}{\check{E}(X_{j\check{\bullet}}^2)} \text{ for }j'\ne j,~\check{U}_{j+1}=\frac{X_{j}X_{j\check{\bullet}}}{\check{E}(X_{j\check{\bullet}}^2)}-1$\\
	 {} & {} &  $\hat{a}, \hat{\bf b}$ are the LSE of model \eqref{eq:model1} with\\
	{} & {} & $(\check{W}_j,\check{U}_{1},\ldots,\check{U}_{p+1})$ replacing $({W}_j,U_{1},\ldots,U_{p+1})$\\
	\hline
	Asymptotics & $\nu =\lim \frac{n}{n+m} $ & $n$ ($m$) the labeled (unlabeled) sample size \\
	{} & ${\boldsymbol \Sigma}_{LSE}$, ${\boldsymbol \Sigma}_{PI}$, ${\boldsymbol \Sigma}_{TI}$ & asymptotic covariance matrix of LSE, PI, TI \\
	{} & {} & as in Theorem \ref{thm:p_dim}\\
	\hline 
	Main results {}&  ${\boldsymbol \Sigma}_{LSE} =  {\boldsymbol \Sigma}_{TI}+{\boldsymbol \Sigma}_{\rm diff}$ & ${\boldsymbol \Sigma}_{\rm diff}=Cov(\delta {\bf X}_{\bullet}-\tilde{\boldsymbol \delta})$\\
	{} & ${\boldsymbol \Sigma}_{PI}= {\boldsymbol \Sigma}_{TI}+\nu{\boldsymbol \Sigma}_{\rm diff}$ & {} \\
	\hline
 		\end{tabular}}
\label{tab:not2}
\end{center}
\end{table}

Let $\vec{\bf X}_{-j} = (1,X_1,\ldots,X_{j-1},X_{j+1},\ldots,X_p)^T$ and let 
\begin{equation}\label{eq:proj}
{\boldsymbol \beta}_{-j \bullet } = \arg \min_{\tilde {\boldsymbol \beta}} E \left(X_j- \tilde {\boldsymbol \beta}^T \vec{\bf X}_{-j}\right)^2 = \left\{E\left( \vec{\bf X}_{-j} \vec{\bf X}_{-j}^T \right) \right\}^{-1} E\left( \vec{\bf X}_{-j} X_j \right).
\end{equation}  
Now, define $X_{j\bullet} = X_j -  \vec{\bf X}_{-j}^T {\boldsymbol \beta}_{-j \bullet }$. 
Each $\beta_j$ can be written in the one dimensional form
$
\beta_j = \frac{E (Y X_{j\bullet} )}{E(X_{j\bullet}^2)}.
$

The standard LSE can be viewed in a similar manner. Let ${\bf Y}=(Y^{(1)},\ldots,Y^{(n)})^T,~{\bf X}_j=(X_j^{(1)},\ldots,X_j^{(n)})^T$ and let ${\bf X}_{-j}= ( {\bf 1}, {\bf X}_1, \ldots, {\bf X}_{j-1},{\bf X}_{j+1},\ldots, {\bf X}_{p} )$. Define also $\hat{ \boldsymbol \beta}_{-j \bar \bullet }= \left\{ {\bf X}_{-j}^T {\bf X}_{-j} \right\}^{-1} {\bf X}_{-j}^T {\bf X}_{j}$ and ${\bf X}_{j \bar{\bullet}} = {\bf X}_j - {\bf X}_{-j} \hat{ \boldsymbol \beta}_{-j \bar \bullet }$; then,
$\{{\hat {\boldsymbol \beta}}_{LSE}\}_j = \frac{\langle \bf{Y},{\bf X}_{j \bar{\bullet}}   \rangle}{|| {\bf X}_{j \bar{\bullet}} ||^2}$. {Recall that $\bar{\cdot}$ denotes mean over the labeled sample; thus, $\bar{\bullet}$ denotes adjustments over the labeled sample, whereas $\bullet$ denotes adjustments over the population.}

The total information estimator is the intercept estimator of the regression model in the population obtained by multiplying \eqref{eq:model} by $\frac{X_{j\bullet}}{E(X_{j\bullet}^2)}$, that is,
\begin{equation}\label{eq:model1}
W_j = \beta_j +a U_1 + \sum_{j'=1}^p b_{j'} U_{j'+1}  + \tilde{\delta}_j,
\end{equation}
where
\begin{equation}\label{eq:model1_def}
W_j = \frac{YX_{j\bullet}}{E(X_{j\bullet}^2)},~U_1 =\frac{X_{j\bullet}}{E(X_{j\bullet}^2)},~U_{j'+1}=\frac{X_{j'}X_{j\bullet}}{E(X_{j\bullet}^2)} \text{ for }j'\ne j,\text{ and } U_{j+1}=\frac{X_{j}X_{j\bullet}}{E(X_{j\bullet}^2)}-1
\end{equation}
are $Y,1,X_1,\cdots,X_p$ multiplied by $\frac{X_{j\bullet}}{E(X_{j\bullet}^2)}$ and $a, {\bf b}$ are $\theta_2,\ldots,\theta_{p+1}$ defined by \eqref{eq:best_lin} for $W=W_j$ and ${\bf U}=(U_1,\ldots,U_{p+1})$.
Since $E(U_1)=\cdots=E(U_{p+1})=0$, then the intercept is $E(W_j)=\beta_j$.
{ In setting $E(W_j)=\beta_j$ and also the expectations of the $U$'s to be zero, we exploited the knowledge of the moments of ${\bf X}$, yielding higher efficiency as shown below.}  
The remainder term in \eqref{eq:model1}, $\tilde{\delta}_j=W_j - \beta_j - a U_1- \sum_{j'=1}^p b_j' U_{j'+1}$ is orthogonal to $(U_1,\ldots,U_{p+1})$.
Notice that the vector $(U_1,\ldots,U_{p+1})$ depends on $j$ but this is suppressed in the notation.

Specifically, we define the total information (TI) estimator to be
\begin{equation*}
\{\hat{\boldsymbol \beta}_{TI}\}_j = \bar{W}_j -\hat{a} \bar{U}_1 -\sum_{j'=1}^p\hat{b}_{j'} \bar{U}_{j'+1},
\end{equation*}
where $\hat{a}, \hat{\bf b}$ are the regular LSE of \eqref{eq:model1}, and $\bar{\cdot}$ denotes the mean over the labeled sample with $n$ observations.

To define the partial information estimator let
\begin{equation}\label{eq:beta_minus_j}
{\boldsymbol \beta}_{-j \check{\bullet} } =  \left\{\check{E}\left( \vec{\bf X}_{-j} \vec{\bf X}_{-j}^T \right) \right\}^{-1} \check{E}\left( \vec{\bf X}_{-j} X_j \right),
\end{equation}  
where $\check{E}$ is the the empirical mean based on the full ${\bf X}$ sample of size $n+m$. Now define $X_{j \check{\bullet}}= X_j- \left\{ \vec{\bf X}_{-j}\right\}^T {\boldsymbol \beta}_{-j \check{\bullet} }$ and
\begin{equation*}
\check{W}_j = \frac{YX_{j\check{\bullet}}}{\check{E}(X_{j\check{\bullet}}^2)},~\check{U}_1 =\frac{X_{j\check{\bullet}}}{\check{E}(X_{j\check{\bullet}}^2)},~\check{U}_{j'+1}=\frac{X_{j'}X_{j\check{\bullet}}}{\check{E}(X_{j\check{\bullet}}^2)} \text{ for }j'\ne j,\text{ and }\check{U}_{j+1}=\frac{X_{j}X_{j\check{\bullet}}}{\check{E}(X_{j\check{\bullet}}^2)}.
\end{equation*}
{Here ${\check{\bullet}}$ denotes adjustments over the full $n+m$ sample. }
The partial information (PI) estimator is
\begin{equation*}\label{eq:PI}
\{\hat{\boldsymbol \beta}_{PI}\}_j = \bar{\check{W}}_j -\hat{a} \bar{\check{U}}_1 -\sum_{j'=1}^p\hat{b}_{j'} \bar{\check{U}}_{j'+1},
\end{equation*}
where $\hat{a}, \hat{b}_{j'}$ are the regular LSE of model \eqref{eq:model1} with $(\check{W}_j,\check{U}_{1},\ldots,\check{U}_{p+1})$ replacing $({W}_j,U_{1},\ldots,U_{p+1})$.

The following theorem establishes the asymptotic distribution of the estimates and states conditions under which $\hat{\boldsymbol \beta}_{TI}$ and $\hat{\boldsymbol \beta}_{PI}$, asymptotically dominate $\hat{\boldsymbol \beta}_{LSE}$.
The asymptotic distribution of $\hat{\boldsymbol \beta}_{LSE}$ is already known \citep{White1} and it is presented here so that the comparison to $\hat{\boldsymbol \beta}_{TI}$ and $\hat{\boldsymbol \beta}_{PI}$ can be made.

\begin{theorem} \label{thm:p_dim}
\begin{enumerate}
\item[(i)] Suppose that the vector ${\bf X}_{\bullet}=\left\{\frac{X_{1\bullet}}{E(X_{1 \bullet}^2)},\ldots,\frac{X_{p\bullet}}{E(X_{p \bullet}^2)}\right\}$ is well defined (i.e., the projections in \eqref{eq:proj} exist and $E(X_{j \bullet})^2$ is positive and finite for $j=1,\ldots,p$) and that  $\delta {\bf X}_{\bullet}$ has finite second moments,  then,
\[
\sqrt{n}(\hat{\boldsymbol \beta}_{LSE}-{\boldsymbol \beta}) \tendd N(0,{\boldsymbol \Sigma}_{LSE}),
\]
where 
${\boldsymbol \Sigma}_{LSE}= Cov(\delta {\bf X}_{\bullet})$.
\item[(ii)]  Suppose further that for each $j=1,\ldots,p$, the vector $(W_j,U_1,U_2,\ldots,U_{p+1})$ has finite second moments and that the matrix $E({\bf U}_j {\bf U}_j^T)$ is invertible for ${\bf U}_j= \left( 1, U_1, U_2,\ldots,U_{p+1} \right)^T$;
then,
\[
\sqrt{n}(\hat{\boldsymbol \beta}_{TI}-{\boldsymbol \beta}) \tendd N(0,{\boldsymbol \Sigma}_{TI}),
\]
where ${\boldsymbol \Sigma}_{TI}= Cov(\tilde{\boldsymbol \delta})$ for $\tilde{\boldsymbol \delta}=(\tilde{\delta}_1,\ldots,\tilde{\delta}_p)$. Furthermore, we have that ${\boldsymbol \Sigma}_{LSE} =  {\boldsymbol \Sigma}_{TI}+{\boldsymbol \Sigma}_{\rm diff}$ where 
${\boldsymbol \Sigma}_{\rm diff}=Cov({\bf V})$ for ${\bf V}=\delta {\bf X}_{\bullet}-\tilde{\boldsymbol \delta}$.
\item[(iii)] Suppose further that $\lim \frac{n}{n+m} = \nu$, then, 
\[
\sqrt{n}(\hat{\boldsymbol \beta}_{PI}-{\boldsymbol \beta}) \tendd N(0,{\boldsymbol \Sigma_{PI}}),
\]
where ${\boldsymbol \Sigma}_{PI}= {\boldsymbol \Sigma}_{TI}+\nu{\boldsymbol \Sigma}_{\rm diff}$.
\end{enumerate}
Therefore, if ${\boldsymbol \Sigma}_{\rm diff}$ is not the zero matrix, then ${\boldsymbol \Sigma}_{LSE}-{\boldsymbol \Sigma}_{TI}$ is positive definite, and if further $\nu<1$, then ${\boldsymbol \Sigma}_{LSE}-{\boldsymbol \Sigma}_{PI}$ is also positive definite.
\end{theorem}
\begin{corollary} \label{cor:p_dim}
Theorem \ref{thm:p_dim} and Proposition \ref{prop1} imply that $R(\hat{\alpha}_{LSE},\hat{\boldsymbol \beta}_{LSE})- R(\hat{\alpha}_{PI},\hat{\boldsymbol \beta}_{PI}) = \text{Trace}\{(1-\nu){\boldsymbol \Sigma}_{\rm diff} {\bf M}\}$, where ${\bf M}=Cov({\bf X})$.
\end{corollary}

In short, Theorem \ref{thm:p_dim} states that if the regression models \eqref{eq:model} and \eqref{eq:model1} are well defined and the residuals have finite second moments, then ${\boldsymbol \Sigma}_{LSE}-{\boldsymbol \Sigma}_{TI}$ is positive semi-definite and it is strictly positive unless ${\bf V} \equiv 0$. If further, the unlabeled sample size is not negligible, i.e. $\lim m/n > 0$, then ${\boldsymbol \Sigma}_{LSE}-{\boldsymbol \Sigma}_{PI}$ is also positive definite.

As in the one-dimensional case, $\hat{\boldsymbol \beta}_{TI}$ and $\hat{\boldsymbol \beta}_{PI}$ improve over $\hat{\boldsymbol \beta}_{LSE}$ only when ${\boldsymbol \Sigma}_{\rm diff}$ is not the zero matrix or equivalently that ${\bf V}$ is not zero. For each $j$, $Var(V_j)$ measures the difference between the original model \eqref{eq:model} and the intercept model \eqref{eq:model1}. When 
\begin{equation}\label{eq:cond_linear}
E(X_{j\bullet}^2 \delta)=E(X_{j\bullet}^2 X_1 \delta)=\cdots=E(X_{j\bullet}^2 X_p \delta)=0 
\end{equation}
then $\delta X_{j\bullet} \equiv \delta_j$ and the two models are essentially the same, in which case $V_j \equiv 0$. Otherwise, if for some $j$ \eqref{eq:cond_linear} does not hold, then $\hat{\boldsymbol \beta}_{PI}$ improves over $\hat{\boldsymbol \beta}_{LSE}$ when $\nu <1$; see Appendix \ref{app:F}. {In the one-dimensional case, \eqref{eq:cond_linear} is equivalent to $E(X^2 \delta)=E(X^3 \delta)=0$, i.e., that $X^2$ and $X^3$ are uncorrelated with $\delta$. Generally, \eqref{eq:cond_linear} implies that certain linear combinations of $E(X_j X_{j'} \delta)$ and $E(X_{j} X_{j'} X_{j''} \delta)$ are 0. Thus, roughly speaking, {when $E[Y | {\bf X}]$ is non-linear and the non-linear part is correlated with several second or third moments of ${\bf X}$, then} we expect improvement of $\hat{\boldsymbol \beta}_{PI}$.}

\subsection{New methodology (summary)}

We now provide a step-by-step description of the new methodology for estimation of ${\boldsymbol \beta}$.

For each $j=1,\ldots,p$:
\begin{enumerate}
\item (Adjust the regressors) Let ${\boldsymbol \beta}_{-j \check{\bullet} }$ as defined in \eqref{eq:beta_minus_j} and define $X_{j \check{\bullet}}^{(i)}= X_j^{(i)}- \left\{ \vec{\bf X}_{-j}^{(i)}\right\}^T {\boldsymbol \beta}_{-j \check{\bullet} }$ for $i=1,\ldots,n$.
\item (Define the intercept model) Define for $i=1,\ldots,n$
\begin{equation*}
\check{W}_j ^{(i)} = \frac{Y^{(i)}X_{j\check{\bullet}}^{(i)}}{\check{E}(X_{j\check{\bullet}}^2)},~\check{U}_1^{(i)} =\frac{X_{j\check{\bullet}}^{(i)}}{\check{E}(X_{j\check{\bullet}}^2)},~\check{U}_{j'+1}^{(i)}=\frac{X_{j'}^{(i)}X_{j\check{\bullet}}^{(i)}}{\check{E}(X_{j\check{\bullet}}^2)} \text{ for }j'\ne j,\text{ and }
\check{U}_{j+1}^{(i)}=\frac{X_{j}^{(i)}X_{j\check{\bullet}}^{(i)}}{\check{E}(X_{j\check{\bullet}}^2)}-1.
\end{equation*}
\item (Define the intercept estimator) The partial information (PI) estimator is
\begin{equation*}\label{eq:PI}
\{\hat{\boldsymbol \beta}_{PI}\}_j = \bar{\check{W}}_j -\hat{a} \bar{\check{U}}_1 -\sum_{j'=1}^p\hat{b}_{j'} \bar{\check{U}}_{j'+1},
\end{equation*}
where $\hat{a}, \hat{b}_{j'}$ are the regular LSE of the regression model
\begin{equation*}
W_j^{(i)} = \beta_j +a \check{U}_1^{(i)} + \sum_{j'=1}^p b_{j'} \check{U}_{j'+1}^{(i)} + \tilde{\delta}_j^{(i)},\text{for }i=1,\ldots,n.
\end{equation*}
\end{enumerate}
{The methodology is built from standard least squares procedures. An R code that implements the algorithm and also computes estimates of the variance as in Section \ref{sec:est_var} below, is available at the homepage of the first author.}

\section{{Simplicity and semi-parametric efficiency}} \label{sec:simple}

{This section describes an optimality property of $\hat{\beta}_{TI}$ and discusses its relation to the estimator of \citet{Chakrabortty} and to the notion of semi-parametric efficiency.}   

\subsection{Adding terms to the linear model}

Our methodology is based on replacing the original linear model \eqref{eq:model1_1} with the intercept model \eqref{eq:modelXY}. While the original model has two covariates $(1,X)$ in the one-dimensional case, the intercept model has three covariates $(1,U_1,U_2)$. Therefore, besides of replacing the model, our method adds one covariate to the model. In Section \ref{sec:simp} below we show that this addition is necessary and sufficient for improving the LSE asymptotically. Here we argue that adding terms in the original linear model does not always lead to improvement over the LSE. That is, adding orthogonal terms to the regression model can sometimes harm the estimation of the parameters. 

	
To observe the latter fact, consider the one dimensional case (Model \eqref{eq:model1_1}) and the total information scenario where $E(X)=0$ and $E(X^2)=1$. One could modify the model to
	\begin{equation} \label{eq:modify_model}
	Y= \alpha + \beta X+ f(X) + \tilde{\delta},
	\end{equation}
	where $f(X)$ is some function of $X$ that satisfies $E(f(X))=0$ and $E(X f(X))=0$ (otherwise $\alpha$ and $\beta$ might change). In this case, the asymptotic variance for the LSE of $\beta$ in the original model \eqref{eq:model1_1} is $E(X^2 \delta^2)$ and in the modified model \eqref{eq:modify_model} it is $E(X^2 \tilde{\delta}^2)$. For a given $f(X)$, it is  possible, under mild conditions, to find ${\delta}$ such that the modified estimator has larger asymptotic variance. To construct one, pick an arbitrary residual $\tilde{\delta}=\tilde{\delta}_0$ such that $E (X^2 f(X) \tilde{\delta}_0)$ is nonzero (assuming such $\tilde{\delta}_0$ exists). Then in the original model \eqref{eq:model1_1}, consider the family $\delta=\delta_c=f(X)+c \tilde{\delta}_0$ for variable $c$. The difference of the asymptotic variances between the new LSE and the old LSE is 
	\[
	E(X^2 c^2 \tilde{\delta}_0^2)-E(X^2 \delta_c^2)=-E(X^2 f^2(X))-2c E (X^2 f(X) \tilde{\delta}_0).
	\]
	The variable $c$ can be chosen to make this difference positive, yielding a model in which the ordinary least squares estimation has been harmed.

\subsection{{Semi-parametric efficiency}}

{To define a semi-parametric efficient estimator in this context, we rewrite model \eqref{eq:model} as
\begin{equation}\label{eq:gmodel}
Y=\alpha+{\bs \beta}^T {\bf X} + \underbrace{\eta + \ve}_{=\delta},
\end{equation}
where $\eta=\eta({\bf X}) = E(Y|{\bf X})-\alpha - {\bs \beta}^T {\bf X}$ and $\ve=Y-\alpha-{\bs \beta}^T {\bf X} - \eta$ and assume that the marginal distribution of ${\bf X}$ is known (total information). 
In this model, the semi-parametric efficiency bound for the variance is $\left\{E(\vec{\bf X} \vec{\bf X}^T)\right\}^{-1}E(\vec{\bf X} \vec{\bf X}^T \ve^2) \left\{E(\vec{\bf X} \vec{\bf X}^T)\right\}^{-1}$ (see \citet{Chakrabortty} for more details). } 

{Here is a general description of the approach of Chakrabortty and Cai, which is based on the idea of imputation of the $Y$'s in the unlabeled sample. Specifically, when $Y_i$ is missing, define $\hat{Y}_i = \hat{m}({\bf X}_i)$ for $i=n+1,\ldots,n+m$, where $\hat{m}$ is an estimate of $m({\bf X})=E(Y|{\bf X})=\alpha+{\bs \beta}^T{\bf X} +\eta(\bf X)$ based on the labeled sample. They suggest to use the LSE of the unlabeled sample when using the imputed $Y$'s, i.e., based on $({\bf X}_n,\hat{Y}_n),\ldots,({\bf X}_{n+m},\hat{Y}_{n+m})$. When $\hat{m}$ converges to $m$, the resulting estimator is semi-parametric efficient. However, in some situations, e.g., when ${m}(\bf X)$ has a complicated parametric form and $n$ is not very large, it is difficult to efficiently estimate it with the dataset in hand. Below we review the pros and cons of this approach with respect to ours. }


\subsection{{Simplicity}} \label{sec:simp}

{We now discuss the relation of our approach to semi-parametric efficiency and introduce an optimality property of our approach, which we call ``simplicity''. 
In order to avoid unnecessary technicalities, we consider here only the one-dimensional case under full information. Explicitly, let $(X,Y) \sim G$ where
\begin{equation} \label{eq:cond}
Var_G(X) \in (0,\infty),~E_G(XY)^2 < \infty \text{ and } E_G(X)=E_G(X^2)-1=0; 
\end{equation}
the latter assumption means that the first and second moments of $X$ are known. The sub-script $G$ is used to emphasize that the expectation is taken over $G$ as we will consider below several such $G$'s. In this context, the semi-parametric efficiency bound for the variance is $E_G(\ve^2 X^2)$. } 

{We argue that $\hat{ \beta}_{TI}$ is the simplest (to be defined shortly) estimator within the framework of intercept model that dominates the LSE, and we also describe how to construct an estimator that achieves the semi-parametric efficiency bound, which is different from the estimator of Chakrabortty and Cai. To this end, consider a basis of functions $\{q_k(X)\}_{k=1}^\infty$ that satisfies  
\begin{multline} \label{eq:cond1}
\text{for all }K\in \mathbb{N}:~ E_G[q_1(X)]=\cdots=E_G[q_K(X)]=0,\\
E_G\left[\{1,q_1(X),\ldots,q_K(X)\}^T \{1,q_1(X),\ldots,q_K(X)\}\right]\text{ is invertible,}\\
\text{and }\{XY,q_1(X),\ldots,q_K(X)\}\text{ has finite second moments under }G. 
\end{multline}
The assumption that $E_G[q_1(X)]=\cdots=E_G[q_K(X)]=0$ is justified since the marginal distribution of $X$ under $G$ is assumed known.  }

{Consider an intercept model based on $q_1(X),\ldots,q_K(X)$
\begin{equation}\label{eq:int_model}
XY=\beta + \sum_{k=1}^K b_k q_k(X) + \tilde{\delta}_K.
\end{equation}
When \eqref{eq:cond1} holds, then $b_1,\ldots,b_k$ in \eqref{eq:int_model} are well defined as the parameters of the best linear predictor. Since  $E_G[q_1(X)]=\cdots=E_G[q_K(X)]=0$, the intercept is $E_G(XY)=\beta$.}   

{Suppose we have $n$ iid observations $(X_1,Y_1),\ldots,(X_n,Y_n)$ from $G$. Let $\hat{\beta}_{LSE}$ be the least squares estimate based on $n$ observations from Model \eqref{eq:model_1} and let $\hat{\beta}_{K}$ be the intercept LSE estimate based on $n$ observations from Model \eqref{eq:int_model}. By the ``Sandwich'' Theorem \eqref{eq:sandV},
\[
\sqrt{n}\left(\hat{\beta}_{LSE} - \beta \right) \tendd N\left(0, E_G(\delta X)^2\right)\text{ and }\sqrt{n}\left(\hat{\beta}_{K} - \beta \right) \tendd N\left(0, E_G\tilde{\delta}_K^2\right).
\] 
The following theorem states that $\hat{\beta}_{K}$ is better (smaller asymptotic variance) than $\hat{\beta}_{LSE}$ for every distribution $G$, iff $X \in \text{span}\{q_1(X),\ldots,q_K(X)\}$ and $X^2-1 \in \text{span}\{q_1(X),\ldots,q_K(X)\}$.}  

{
	\begin{theorem}\label{thm:min}
	${}$	\\
	(Sufficiency) If $X \in \text{span}\{q_1(X),\ldots,q_K(X)\}$ and $X^2-1 \in \text{span}\{q_1(X),\ldots,q_K(X)\}$, then $\hat{\beta}_{K}$ is better than $\hat{\beta}_{LSE}$  for ``every'' $G$, i.e., $E_G\tilde{\delta}_K^2 \le E_G(\delta X)^2$ for every distribution $G$ that satisfies \eqref{eq:cond} and \eqref{eq:cond1}.\\
	(Necessity) If $X \notin \text{span}\{q_1(X),\ldots,q_K(X)\}$ or $X^2-1 \notin \text{span}\{q_1(X),\ldots,q_K(X)\}$, then there exists distribution $G$ that satisfies \eqref{eq:cond} and \eqref{eq:cond1}, and where $\hat{\beta}_{LSE}$ is better than $\hat{\beta}_{K}$, i.e., $E_G(\delta X)^2 < E_G\tilde{\delta}_K^2$.
\end{theorem}}

{Theorem \ref{thm:min} implies that $\hat{\beta}_{TI}$ is ``simple'' in the sense that it is derived from the minimal intercept model that is better than $\hat{\beta}_{LSE}$ for every $G$ that satisfies \eqref{eq:cond} and \eqref{eq:cond1}. In this sense, $\hat{\beta}_{TI}$ is the simplest estimator that is derived from an intercept model.}

{This ``simplicity'' property allows us to improve over the LSE without explicitly modeling the non-linear part $\eta({\bf X})$. This is one advantage over Chakrabortty and Cai's approach, where the missing $Y$'s are imputed using an estimate of $\eta({\bf X})$. We expect that if $\eta({\bf X})$ has a complicated form and $n$ is not very large, our estimate may be advantageous. On the other hand, a successful approximation of $\eta({\bf X})$, may lead to an  estimator of ${\bs \beta}$ with smaller variance, close to the semi-parametric bound. This is demonstrated and further discussed in Section \ref{sec:comp_sim} below.} 

{The intercept model framework can also be used to construct an estimator that achieves the semi-parametric bound. To this end, assume for the moment that for some $K$   
\begin{equation}\label{eq:span}
\alpha X + \beta X^2 + \eta X \in \text{span}\{q_1(X),\ldots,q_K(X)\}.
\end{equation}
Then Model \eqref{eq:gmodel} implies that $E_G(\tilde{\delta}_K^2)=E_G(\ve^2 X^2)$ and therefore $\hat{\beta}_K$ is semi-parametric efficient. Generally, \eqref{eq:span} is not likely to hold precisely, but if $\{q_k(X)\}_{k=1}^\infty$ is rich enough, then for large enough $K$, \eqref{eq:span} holds true up to some small enough approximation error. In this case, one can increase $K$ slowly with $n$ and then achieve semi-parametric efficiency. This procedure is made precise in Theorem 7 of \citet{Anru}.}

 {To sum up the findings of this section we  conclude that the intercept model framework derives a large class of estimates that are better than $\hat{ \beta}_{LSE}$.  Our suggested estimator $\hat{ \beta}_{TI}$ is minimal within this class and a semi-parametric efficient estimator can be achieved by considering a sequence of increasing models. Compared to Chakrabortty and Cai's estimator it is especially useful when the non-linear part is difficult to model.}

\section{A simulation study} \label{sec:sim}

We compare the performance of the Partial Information (PI) and Total Information (TI) estimators against the Least Squares Estimate (LSE) across a wide range of settings. 

\subsection{Toy example} 

We start by studying the following toy model
\begin{equation} \label{eq:toymodel}
Y=\alpha X^2 + \beta X + \ve = \alpha + \beta X + \underbrace{\alpha(X^2-1) + \ve}_{=\delta},
\end{equation}
where $X$ and $\ve$ are i.i.d N(0,1). Under model \eqref{eq:toymodel}, the linear and non-linear part are determined by $\beta$ and $\alpha$ separately and the linear coefficient $\beta$ does not affect the residual $\delta$. In this case, it easy to calculate the asymptotic variance of the estimates explicitly 
\[
\sigma^2_{LSE}=10\alpha^2+1,~\sigma^2_{TI}=6 \alpha^2 +1,~\sigma^2_{\rm diff}=4 \alpha^2.
\] 
Considering the excess risk, then by Proposition \ref{prop1} the ratio of excess risks of $\hat{\beta}_{PI}$ ($\hat{\beta}_{TI}$, respectively) and $\hat{\beta}_{LSE}$ converges to $\frac{2\alpha^2+1+\sigma^2_{PI}}{2\alpha^2+1+\sigma^2_{LSE}}$ ($\frac{2\alpha^2+1+\sigma^2_{TI}}{2\alpha^2+1+\sigma^2_{LSE}}$, respectively). The limits equal 1 in the linear case ($\alpha=0$) and approach 0 when the non-linear part is dominant ($\alpha \to \infty$).  

Figure \ref{fig:toy} summarizes 10,000 simulations of model \eqref{eq:toymodel} with $\beta=1$ and different values of $\alpha$, $n$; for the PI case we used $m=2n$. 
For all the scenarios we found that $\hat{\beta}_{TI}$, $\hat{\beta}_{PI}$ have smaller excess risk when $n$ is large enough. As $\alpha$ increases,  the departure from linearity is more significant and the ratio of the excess risks is smaller. For small $n$, LSE is superior for all scenarios as it is a simpler estimate. When the model is close to linear ($\alpha=1/4$), the new estimates are better for $n\ge 200$ and for the other values of $\alpha$ the new estimates are better for $n \ge 70$. For small values of $\alpha$ the limiting excess risk is close to the actual risk even when $n\approx 300$ but this does not hold true for larger values of $\alpha$.
In short, we found that the new estimates are better for large $n$, and are much better when the non-linear part is significant.

\begin{figure}[ht!] 
        \subfigure[TI ]{%
            \includegraphics[width=0.5\textwidth]{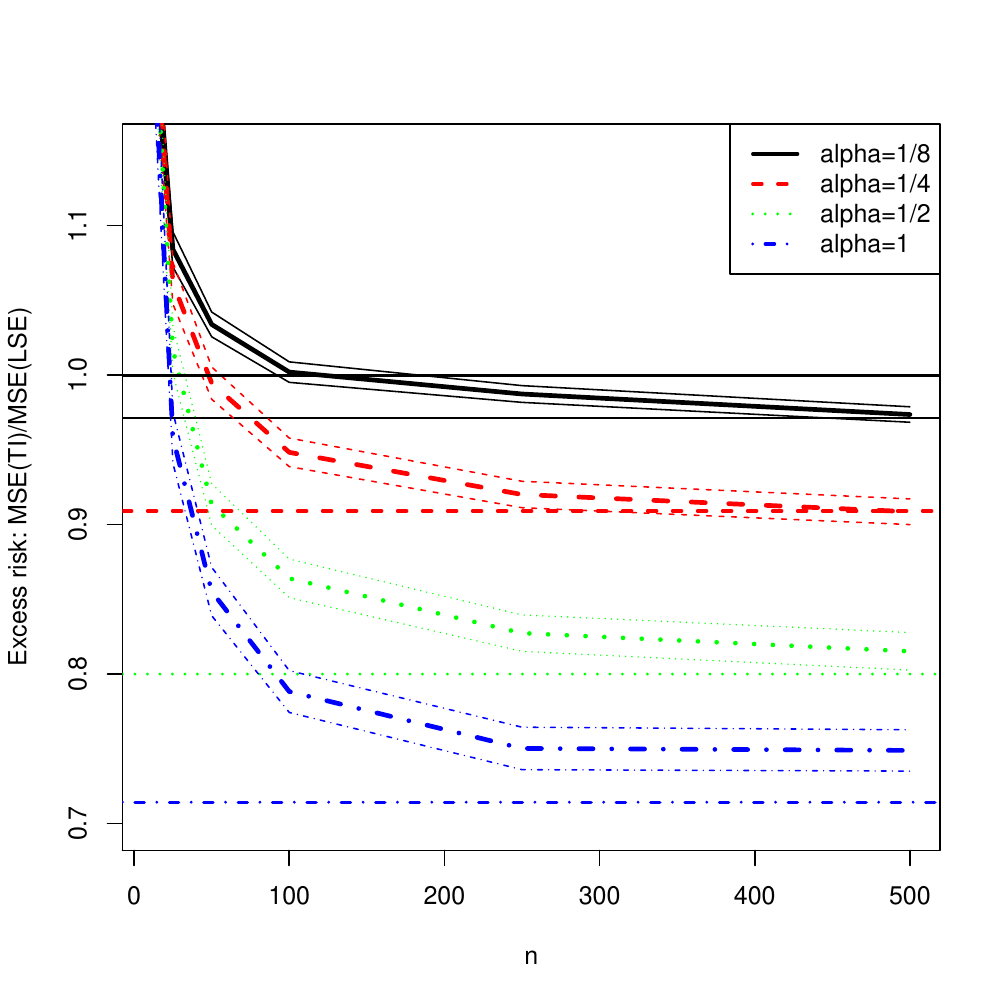}
        }
        \subfigure[PI]{%
           \includegraphics[width=0.5\textwidth]{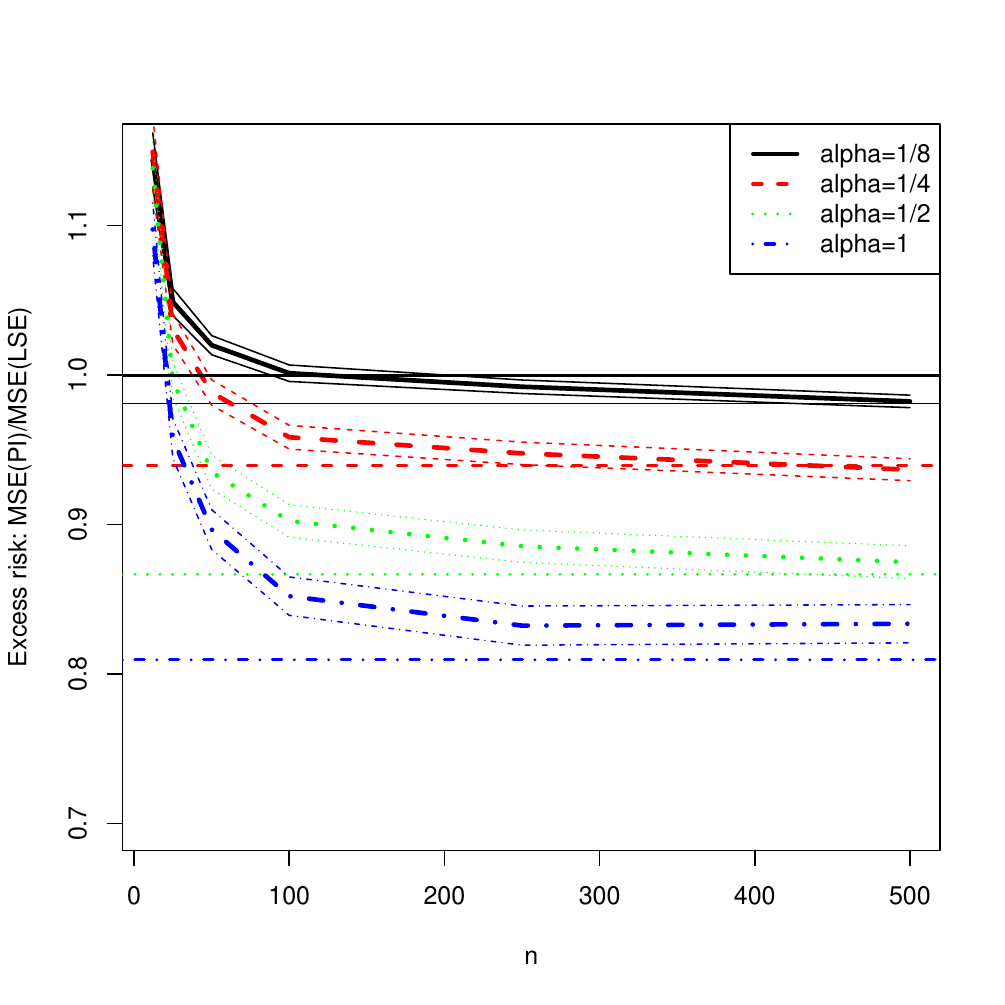}
        }\\    \caption{\footnotesize Estimates of the excess risk ratio for the TI and PI estimates based on 10,000 simulations of model \eqref{eq:toymodel} with $\beta=1$ and different values of $\alpha$, $n$; for the PI case we used $m=2n$. A confidence interval based on two standard deviations is also plotted. The horizontal lines represents the limiting excess risk. \label{fig:toy}}
\end{figure}
 
\subsection{Simulation Parameters}\label{sec:sim_par}

We now compare our new estimates to LSE for a broader variety of scenarios.
Estimates of the Excess Risk Ratio (ERR) for PI and TI were obtained for every combination of the following parameter choices:

\bigskip

\begin{tabular}{ l r }
	$n$ (size of labeled dataset) & 12, 25, 50, 100, 250, 500 \\
	$m$ (size of unlabeled dataset) & $n, 2n$ \\
	$p$ (number of predictors) & 1, 4 \\
	Distribution of $X$ & Gaussian, Lognormal, Exponential, Cubed Gaussian \\
	Errors & Gaussian, $N(0,e^{2\|X\|})$ \\
	$E[Y|X]$ & $X$, $e^X$, $X^3$, $ \sqrt{|X|}$ \\
\end{tabular}

\bigskip

{For $p = 4$, $E[Y|X]$ was the sum over $j$ of the chosen function of each $X_j$}. Fixing each of the 768 parameters settings, many sample datasets were generated, along with a large test dataset (n = 100,000). For each sample dataset, PI, TI, and LSE were fit, (calculating TI using the PI method with $m \geq 500 n$.) Next, an estimate of expected MSE was obtained using the test dataset. Along with the best linear fit for the test dataset (given by least squares regression), these determine the Excess Risk Ratio for PI ($ERR_{PI}$) as follows: 
\[
ERR_{PI} = \frac{MSE_{PI} - MSE_{BLF}}{MSE_{LSE} - MSE_{BLF}},
\]
where \\
$MSE_{PI}$ is the expected MSE of the PI estimator (estimated from mean performance on the test dataset), \\
$MSE_{LSE}$ is the expected MSE of the LSE estimator, and \\
$MSE_{BLF}$ is the MSE of the best linear fit of the test dataset.

An analagous calculation was performed to calculate $ERR_{TI}$. Additional datasets were sampled to improve estimates of $MSE_{PI}$ and $MSE_{LSE}$, until standard errors for $ERR_{PI}$ and $ERR_{TI}$ fell below $1\%$ (as estimated by the delta method), or a maximum of $100,000$ sample datasets was generated.

\subsection{Simulation Results}

Table \ref{tab:prop} provides the proportions of scenarios where the PI estimate yields a statistically significant smaller excess risk than the LSE and the proportion where the opposite holds true, for different $p$ and $n$. Statistical significance is measured by two standard errors. The results demonstrate that the PI estimate outperforms the LSE across a wide range of scenarios for large $n$, and the proportion increases with $n$. When $n=500$, LSE is significantly better only for about 10\% of the scenarios we studied. More detailed comments are given in Appendix \ref{sec:sims}:

\begin{table} [ht!] 
\caption{\footnotesize Proportions of scenarios where $ERR_{PI}$ is significantly smaller than 1 and where it is significantly larger than 1. }
\begin{center}
{
  \begin{tabular}{c||c|| c c c c c c}
    \hline
		 {} &\backslashbox{p}{n} & 12 & 25 & 50 & 100 & 250 & 500\\ 
		\hline
	 \multirow{ 2}{*}{$ERR_{PI}<1$ significantly}& 1& 0.297 & 0.344 & 0.469 & 0.594 & 0.656 & 0.609 \\
	 &4 &0.047 & 0.234 & 0.391 & 0.438 & 0.594 & 0.547 \\
	\hline
	\multirow{ 2}{*}{$ERR_{PI}>1$ significantly}& 1&  0.344 & 0.234 & 0.172 & 0.125 & 0.125 & 0.109 \\
	 &4 &0.578 & 0.391 & 0.234 & 0.188 & 0.188 & 0.125 \\
	 \end{tabular}}
\label{tab:prop}
\end{center}
\end{table}

		
		\subsection{A Comparison to other semi-supervised methods}\label{sec:comp_sim}
		
	{
			This section reports a simulation study that compares our estimator, $\hat{\bs \beta}_{PI}$ against three semi-supervised methods: the EASE method of  \citet{Chakrabortty} (using kernel machine regression as the smoothing method), JT-ENET (joint trained elastic net) of \citet{Ryan} and SSR (Semi-supervised regression) from the package SSR in R, which is based on \citet{Hady}. 
			
			Here the dimension is $p=4$, the number of labeled observations is $n=100$, and the unlabeled sample is larger $m=2,000$. For JT-ENET we used only half of the unlabeled data for computational reasons. The simulations were repeated 1000 times for all methods besides JT-ENET, which requires more computational time and therefore was repeated only 100 times. Five scenarios for the joint distribution of $({\bf X},Y)$ were considered; in all scenarios the distribution of ${\bf X}$ and $\ve$ is iid $N(0,1)$.
		\begin{itemize}
		\item {Scenario 1: $Y=(1,1,1,1){\bf X}+ \ve$}
		\item {Scenario 2: $Y=(1,1,1,1){\bf X}+(1,1,1,1){\bf X}^2 + \ve$ (where ${\bf X}^2=(X_1^2,X_2^2,X_3^2,X_4^2)^T$ and similar notation is used below for other functions).}
		\item {Scenario 3: $Y=(1,1,1,1){\bf X}+(0.3,0.3,0.3,0.3){\bf X}^2 + \ve$.}
		\item {Scenario 4: $Y=(1,1,1,1){\bf X}+(1,1,1,1)\left({\bf X}^3-{\bf X}^2+\exp(\bf X)\right) + \ve$.}
		\item {Scenario 5: $Y=(1,1,1,1){\bf X}+(0.3,0.3,0.3,0.3)\left({\bf X}^3-{\bf X}^2+\exp(\bf X)\right) + \ve$.}
	\end{itemize}
Scenario 1 is linear while the rest are not. In Scenarios 2 and 3 the non-linear part has a squared form; they differ in the level of non-linearity. Scenarios 4 and 5 have a complex non-linear part, which is difficult to estimate given the relatively small number of observations.

\begin{table} [ht!] 
\caption{{\footnotesize Mean (std) of the $n\times MSE$ of $\hat{\bs \beta}_{LSE},\hat{\bs \beta}_{PI}$ (our method),  EASE and JT-ENET.} }
\begin{center}
	{
		\begin{tabular}{c||c| c| c | c|c}
			\hline
			{} & $\hat{\bs \beta}_{LSE}$ & $\hat{\bs \beta}_{PI}$ & EASE & JT-ENET & SSR \\ 
			\hline
			\hline
			Scenario 1 &  1.03 (0.03) & 1.12 (0.03) & 1.08 (0.03) & 1.02 (0.8) & 1.7 (0.04)\\
			Scenario 2 & 15.9 (0.3) & 13.3 (0.3) & 8.3 (0.3) & 24.0 (2.3) & 14.9 (0.4) \\
			Scenario 3 &  2.3 (0.05) & 2.2 (0.05) & 2.1 (0.05) & 2.6 (0.18) & 3.1 (0.07)\\
			Scenario 4 &  96 (3) & 74 (2) & 92 (3) & 89 (6) & 113 (2)\\
			Scenario 5 &  9.6 (0.3) & 7.7 (0.2) & 9.4 (0.3) & 8.9 (0.6) & 12.8 (0.3)\\
	\end{tabular}}
	\label{tab:sim}
\end{center}
\end{table}

Here we report results in terms of the MSE of the estimates of ${\bs \beta}$ rather than excess risk since it is easier to interpret. (The ERR defined in \ref{sec:sim_par} can be recovered asymptotically as the ratio of the MSE of $\beta$ to the MSE of $\beta_{LSE}$, due to Proposition 1). The mean MSE and standard errors of the simulations are given in Table \ref{tab:sim}. 

In the linear case (Scenario 1) JT-ENET performs similar to the least squares, while the rest of the methods are slightly worst. The other scenarios consider misspecified models and it suits better the framework of our method ($\hat{\bs \beta}_{PI}$) and EASE. 
EASE outperformed all other methods in Scenario 2; its MSE is about half than the LSE. In scenario 3 all methods besides SSR are more or less the same since the non-linear part is relatively small. Scenarios 4 and 5 demonstrate the advantage of our method since the non-linear part is a complex function and is difficult to estimate. In these scenarios, $\hat{\bs \beta}_{PI}$ improves over the LSE and also outperforms the other methods. Notice that the relative improvement of $\hat{\bs \beta}_{PI}$ remains the same across Scenarios 4 and 5 (the MSE in Scenario 5 is about 10\% of the MSE of Scenario 4 in all estimates). {Scenario 4 demonstrates a situation where the non-linear term is quite significant and thus the MSE is large for all of the methods considered. However, the MSE of $\hat{\bs \beta}_{PI}$ is about 20\% (or more) smaller than the other methods.}

			

		 In order to investigate the performance of the estimators on real data we use a bike sharing dataset \citep{Fanaee}.  The data was aggregated on two hourly basis and the aim is to model the number of casual (as opposed to registered) bike rentals (per two hours) as a function of corresponding weather and seasonal information. 
		 To predict the number of bike rentals, a linear model was used with ten covariates, three of which are continuous and the rest are binary. A list of the covariates is given in Appendix \ref{sec:append}.  The number of observations is 17379.
		 		  
		We computed the LSE in the full dataset and considered the resulting estimates to be the true values of the vector ${\bs \beta}$. The covariates were standardized so that the MSE's for different covariates are comparable. The dataset was split at random  into two parts, labeled and unlabeled, where in the unlabeled part the $y$ observations (number of bike rentals) are deleted and treated as missing. The labeled part is of size $n$, for $n=60,70,\ldots,210$, which is much smaller than the number of observations in the full dataset. For each sample size $n$, 1000 random splits were computed and for each split the LSE, and the semi-supervised estimates were calculated (for PI, 2000 splits were computed due to its computational efficiency; for JT-ENET only 100 splits were considered and the unlabeled sample size was 1000 for computational reasons). For each sample size we computed the relative improvement of the semi-supervised method compared to LSE. We computed the average relative improvement over the three continuous covariates since for the binary covariates no improvement is possible because there are no non-linearities (see Appendix \ref{sec:binary}). 
		
		The average relative improvement is plotted in Figure \ref{fig:relat_improve} as a function of $n$. The relative improvement of JT-ENET was smaller compared to the other methods and therefore its simulation results are not plotted. When $n=60,70$ the relative improvement is the same, more or less, across all methods although SSR seems to be slightly better. For $n$ between 80 and 180, the relative improvement  for PI is 5\% - 10\%, and  for EASE and SSR it is lower, but when $n$ is larger than 180, EASE has the largest relative improvement. 
		
		To sum up, in this dataset, when $n$ is smaller than 180, the non-linear part is difficult to model and PI is advantageous; otherwise EASE has a larger relative improvement.  In terms of prediction of the missing values $Y$, there is almost no improvement since most of prediction error comes form the variability of $Y$, rather than from estimation error.

		\begin{figure}[ht!]
			\begin{center}
				\includegraphics[width=0.8\textwidth]{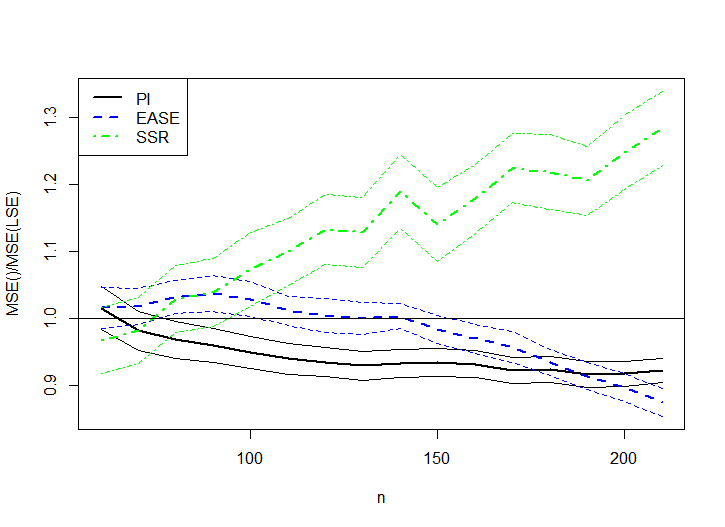}
				\caption{\footnotesize Average relative MSE improvement for PI, EASE and SSR as a function of $n$ in the bike sharing dataset. The thick line denotes the mean improvement over the simulations and the thin line indicates a confidence interval based on two standard deviations.} \label{fig:relat_improve}
			\end{center}
		\end{figure}

		In practice, one does not know the complexity or severity of the nonlinear part relative to $n$, and thus it is unclear which estimator is best. Therefore, it might be advisable to try different methods. Since each method estimates its standard error, one can try to evaluate the proper estimator for the specific dataset in hand.		}


\section{Estimating the standard errors of the estimates} \label{sec:est_se}

In this section we describe a bootstrap estimator of the covariance matrix ${\bs \Sigma}_{PI}$. Other estimates are considered in Appendix \ref{sec:app_est}, where they are compared to the suggested one through simulations.

In the proof of Theorem \ref{thm:p_dim} we showed that for $j=1,\ldots,p$,
\[
\{\hat{\boldsymbol \beta}_{LSE}\}_j=\beta_j+\frac{\frac{1}{n} \sum_{i=1}^n { X}_{j {\bullet}}^{(i)} \delta^{(i)} }{E(X_{j\bullet}^2)}+o_p(1/\sqrt{n})\text{ and }\{\hat{\boldsymbol \beta}_{PI}\}_j = \beta_j+\frac{1}{n}\sum_{i=1}^n \tilde{\delta}^{(i)}+\frac{1}{n+m} \sum_{i=1}^{n+m} V_j^{(i)}+o_p(1/\sqrt{n}).
\]
Hence, the definition of the vector ${\bf V}$ implies that
\begin{multline*}
\hat{\boldsymbol \beta}_{LSE}-\hat{\boldsymbol \beta}_{PI}=\frac{1}{n}\sum_{i=1}^n -{\bf V}^{(i)} +\frac{1}{n+m}\sum_{i=1}^{n+m}{\bf V}^{(i)}+o_p(1/\sqrt{n})\\
= \frac{1}{n}\sum_{i=1}^n {\bf V}^{(i)}\left(\frac{n}{n+m}-1\right)+\frac{1}{n+m}\sum_{i=n+1}^{n+m}{\bf V}^{(i)}+o_p(1/\sqrt{n}). 
\end{multline*}
Therefore, dropping the smaller order terms we obtain
\[
n Cov(\hat{\boldsymbol \beta}_{LSE}-\hat{\boldsymbol \beta}_{PI}) \approx Cov({\bf V})\left(1-\frac{n}{n+m}\right)^2+\frac{nm}{(n+m)^2}Cov({\bf V})
=Cov({\bf V})\left(1-\frac{n}{n+m}\right),
\]
which implies that,
\[
{\boldsymbol \Sigma}_{LSE}-{\boldsymbol \Sigma}_{PI} \approx Cov({\bf V})\left(1-\nu \right).
\]
Thus, a a consistent estimator of the difference ${\boldsymbol \Sigma}_{LSE}-{\boldsymbol \Sigma}_{PI}=:{\boldsymbol \Delta}$ can be obtained by bootstrapping $n Cov ( {\hat {\boldsymbol \beta}}_{PI}- {\hat {\boldsymbol \beta}}_{LSE})$. We denote this estimator by $\widehat{\boldsymbol \Delta}$. 
A consistent estimate to ${{\bs \Sigma}}_{PI}$, which we call variance bootstrap estimator (VBS), is
\[
\hat{{\bs \Sigma}}_{PI,VBS} = \hat{{\bs \Sigma}}_{LSE,BS}- \widehat{\boldsymbol \Delta}, 
\]
where $\hat{{\bs \Sigma}}_{LSE,BS}$ is the standard bootstrap estimator of $\hat{\bs \beta}_{LSE}$ (see details in Appendix \ref{sec:est_var}).
The VBS estimator has the desired property that it is always smaller than the estimated variance of the LSE (in the sense that the difference is positive definite), in a similar fashion to the asymptotic variance.

\section{Inferring homeless population} \label{sec:inferring}

We now consider the Los Angeles homeless dataset \citep{Berk}. Our purpose here is not to carefully analyze this dataset but rather to demonstrate our new method and to compare it to the standard LSE. For our analysis, we used the same covariates as in \citet{Berk}; see this reference for more details about the data.

There are about 2000 census tracts in the Los Angeles county, and a sample was conducted in order to infer the homeless population. The sample consisted of two parts. First, tracts believed to have large numbers of homeless people were pre-selected and visited; there are about 240 such tracts. Second, a sample of about 260 tracts was chosen from the remainder by a random sampling technique and the homeless population was counted, leaving about 1500 tracts to be imputed. 



Leaving aside the preselected tracts this is almost exactly our semi-supervised setting among the remaining population. We have $n=261$ and $m=1536$. The major difference from the setting in our introduction here is that sampling from the remaining population is without replacement and hence is not an i.i.d random sample. Since we have not discussed this type of sampling in detail, we will ignore this issue in the subsequent discussion which is for illustrative purposes. (The difference in sampling schemes does not affect the values of the various estimates, but does impact their standard errors.)

We compare the LSE of the random sample and the new estimate (PI) where the supervised part is the random sample (of size 261) and the unsupervised part consists of the 1536 tracts to be imputed.  The resulting estimates, as well as their standard errors (SE) are given in Table \ref{tab:dataex}. The SE of the LSE was computed using pairs bootstrap ($N_{BS}=10^5$), and the SE of PI was computed using the above variance bootstrap method (the SE of the intercept of PI was estimated by pairs bootstrap). 

We find that the SE's of PI are smaller by 15\%-40\%  than LSE in the predictors MedianInc, PercVacant, PercCommercial, PercIndustrial. These predictors also demonstrated a discrepancy between the SE from linear models to bootstrap SE as reported in \citet{Buja} Section 2, indicating non-linearity in these predictors (a discrepancy was observed also in the predictor PercResidential, but the improvement in this predictor is relatively small). In summary, we found that the additional information on the distribution of the predictors in the unsupervised data provides estimates that are more accurate (smaller variance).

Consider the prediction problem of estimating the homeless population for the 1500 tracts. An estimate of the excess risk ratio, using \eqref{eq:hat_ERR}, is $\widehat{ERR}_{PI}=0.904$, {i.e., improvement of about 10\%. An estimate of the differences of prediction errors as in \eqref{eq:pred_error1} and \eqref{eq:pred_error2} is
$\text{Trace}(\hat{\boldsymbol \Delta} \hat{\bf M})/n=15.2$, while an estimate of the prediction error of the mean is
\[
\frac{1}{m}\sum_{i=1}^m \vec{\bf X}_{n+i}^T \hat{ \boldsymbol \Sigma}_{LSE} \vec{\bf X}_{n+i} /n =42.8.
\]
Since 15.2/42.8 =0.355 the improvement in prediction risk of the means is about 36\%. When considering the prediction of $Y$ itself as in \eqref{eq:pred_error2} then the improvement is much smaller, about 1.5\%.  The reason for this difference is that the predictors we considered are all quite weak, and hence most of prediction error comes from the variability in the distribution of $Y$ given ${\bf X}$ rather than the variability in estimating the $\beta$'s.}

\begin{table} [ht!]
\caption{The estimates $\hat{\boldsymbol \beta}_{LSE}$, $\hat{\boldsymbol \beta}_{PI}$ and their standard errors. }
\begin{center}
{
  \begin{tabular}{c|| c c|c c|c}
    {}& $\hat{\beta}_{LSE}$ &  $\hat{SE}_{LSE}$ & $\hat{\beta}_{PI}$ & $\hat{SE}_{PI}$ & $\frac{{SE}_{LSE}}{{SE}_{PI}}$\\
     \hline
Intercept             & 13.138 & 11.822 & 14.758 & 11.184 & 0.927\\
MedianInc (\$K) 			&-0.065  & 0.056 &-0.080  & 0.044 & 1.284 \\
PercVacant             & 1.449  & 0.707  & 1.583  & 0.514 & 1.374 \\
PercMinority           & 0.060  & 0.105  & 0.063  & 0.099 & 1.058 \\
PercResidential           &-0.070 & 0.100 & -0.078 & 0.094 & 1.072 \\
PercCommercial            & 0.446  & 0.354 & 0.335 & 0.300 & 1.179 \\
PercIndustrial            & 0.101 & 0.188 & 0.202 & 0.143 & 1.317 
\end{tabular}}
\label{tab:dataex}
\end{center}
\end{table}

\section{Discussion} \label{sec:dis}

In this work we showed that ${\boldsymbol \beta}$ can be better estimated when additional information on the predictors ${\bf X}$ is available. The key idea is to replace the regression model 
\begin{equation}\label{eq:model_}
Y= \alpha + \beta_1 X_1+ \cdots +\beta_p X_p + \delta,
\end{equation}
with $p$ regressions of the form
\begin{equation}\label{eq:model1_}
W_j = \beta_j +a U_1 + \sum_{j'=1}^p b_{j'} U_{j'+1} + \tilde{\delta}_j,
\end{equation}
for $j=1,\ldots,p$, where $W_j,U_1,\ldots,U_{p+1}$ is defined in \eqref{eq:model1_def}.
Regression \eqref{eq:model1_} can be used only when the first and second moments of ${\bf X}$ are known. 
We showed that the intercept estimator of \eqref{eq:model1_} is a better estimate of $\beta_j$, in terms of smaller asymptotic variance, than the standard LSE of \eqref{eq:model_}. Furthermore, we showed that even if the second moments of ${\bf X}$ are not known exactly then improvement can be made over the standard LSE provided that the unlabeled sample is non-negligible with respect to the labeled sample, i.e., $m/n$ is bounded away from 0.
 
From a practical viewpoint, in the two datasets we analyzed (in Sections \ref{sec:comp_sim}  and \ref{sec:inferring}), we found that our method improves significantly the estimation of the parameters with respect to LSE but for the purpose of prediction the improvement is small. The reason is that most of the prediction error comes from the variability of $Y$, rather than from estimation error.   

{The improvement of $\hat{\bs \beta}_{PI}$ over $\hat{\bs \beta}_{LSE}$ occurs only when the linear model is misspecified. While there are good reasons to use the linear model even when it is wrong, it should be used with caution (see Section \ref{sec:assumption}). When the primary goal is prediction, then one should compare the linear model to non-linear alternatives using cross-validation, say, and choose the best model. {In this context, sample size is also an important consideration as explained in the next paragraph.} On the other hand, if the primary interest is in understanding the association between the covariates and response, the linear model might be a good choice. When unlabeled data is available the suggested methodology can be helpful in estimating the coefficients of the linear model. } 

{For a specific example, consider the toy model \eqref{eq:toymodel}. Using Proposition \ref{prop1}, it can be shown that the prediction error for the misspecified linear model (with respect to the vector $(1,X)$) is approximately $2\alpha^2+1 +\frac{12\alpha^2+2}{n}$, while the prediction error of the non-linear model (which is a linear model with respect to the vector $(1,X,X^2)$) is $1+\frac{3}{n}$. Notice that the non-linear model here is simple and requires an addition of only one quadratic term, whereas generally non-linear models or non-parametric models are more involved. It follows that if $|\alpha|<1/\sqrt{12}\approx 0.29$, then the linear model has smaller prediction error when $n< \frac{1-12\alpha^2}{2\alpha^2}$. If $|\alpha|\ge 1/\sqrt{12}$ the non-linear model is always better. Thus, for $n_0=\max\{\frac{1-12\alpha^2}{2\alpha^2},0\}$, we have that when $n < n_0$ the linear model has smaller prediction error. The above $n_0$ depends on the amount of non-linearity and when this amount is large $n_0=0$ and the linear model is never better. We expect that this phenomenon occurs more generally, that is, for a typical joint distribution of $({\bf X},Y)$ the linear model is better for a sample size smaller than a certain critical value, which depends on the amount of non-linearity in the conditional expectation $E(Y|{\bf X})$. }

Here we considered the classical framework where $p$ is fixed and $n \to \infty$. \citet{Anru} study the semi-supervised mean estimation problem where $p,n$ both go to infinity but $p$ grows relatively slowly, i.e., $p=O(\sqrt{n})$. The findings of \citet{Anru} indicate that our results can be extended to the latter case. However, our method obviously requires that $p << n$ and for high-dimensional regression models where $p>n$, a different approach is needed. High-dimensional missspecified regression models are the topic of \citet{Buhlman}.

In this work we did not study variable selection and considered the covariates as given. However, improved asymptotic performance may be achieved by including functions of the existing covariates as additional new covariates. Section 3 of \citet{Anru} discusses how this can be done for the problem of semi-supervised estimation of a mean. On the other hand, when the labeled sample, $n$, is small, then it may be better not to include all covariates. We hope to study the problem of variable selection in this context in a future research.

A possible extension of these results is to generalized linear models, where the label variable ($Y$) takes a small number of values and the regression model reduces to a classification problem. Our hope is that our method can be extended to these cases and improvement can be made over naive classifiers that consider only the labeled data. 

In summary, in this work we considered the framework where the best linear predictor is of interest even if $E(Y | {\bf X})$ is not linear.  Under this framework, additional information on the distribution of ${\bf X}$ is useful to construct better estimates than the standard estimates. We believe that this framework is of practical importance and also can lead to interesting statistical issues, which are yet to be studied.         
  
\section*{Acknowledgment}
{
We thank Abhishek Chakrabortty and Mark Culp for helping us with the R code of their methods. We are also grateful to the referees for useful comments that greatly improved the paper. The work of Lawrence D. Brown and Linda Zhao was supported in part by NSF grant DMS-1512084. The work of Michael Sklar was supported in part by NSF grant NSF DMS RTG 1501767  }

\clearpage

\appendix

\noindent{\Large \bf Appendix}

\small

\section{Proofs} \label{sec:proofs}

\subsection*{Proof of Lemma \ref{lem1}}
\begin{lemma}\label{lem1}
We have that
\begin{equation}\label{eq:excess_risk_}
 E L(\tilde{\alpha},\tilde{\boldsymbol \beta})=  E\Big[  \bar{Y} - \sum_{j=1}^p \tilde{\beta}_j\{ {\bar X}_j - E(X_j)\} -E(Y)\Big]^2 +  E\left\{ ({\boldsymbol \beta}- \tilde{\boldsymbol \beta})^T {\bf M} ({\boldsymbol \beta}- \tilde{\boldsymbol \beta}) \right \}.
\end{equation}
\end{lemma}
\noindent{\bf Proof.} \\
First,
\begin{multline*}
R(\tilde{\alpha},\tilde{\boldsymbol \beta})= E\left(Y^* - \tilde{\alpha} - \tilde{\boldsymbol \beta}^T {\bf X}^*\right)^2- E\left(Y^* - {\alpha} - {\boldsymbol \beta}^T {\bf X}^* \right)^2 \\
= E\left\{ (\alpha -\tilde{\alpha}) +   ({\boldsymbol \beta}- \tilde{\boldsymbol \beta})^T{\bf X}^* \right\}\left\{ 2 Y^* - (\alpha +\hat{\alpha}) - 
({\boldsymbol \beta}+ \tilde{\boldsymbol \beta})^T{\bf X}^* \right\}\\
= E\left\{ (\alpha -\tilde{\alpha}) +   ({\boldsymbol \beta}- \tilde{\boldsymbol \beta})^T{\bf X}^* \right\}^2 + 2E\left\{ (\alpha -\hat{\alpha}) +   ({\boldsymbol \beta}- \tilde{\boldsymbol \beta})^T{\bf X}^* \right\}\left\{ Y^* - \alpha  - 
{\boldsymbol \beta}^T{\bf X}^* \right\}  \\
=E\left\{ (\alpha -\tilde{\alpha}) +   ({\boldsymbol \beta}- \tilde{\boldsymbol \beta})^T{\bf X}^* \right\}^2,
\end{multline*}
where the last equality holds true since $Y^* - \alpha  - {\boldsymbol \beta}^T{\bf X}^*$ is orthogonal to $X^*,1$. We can further write
\begin{multline} \label{eq:term1_}
E\left\{ (\alpha -\tilde{\alpha}) +   ({\boldsymbol \beta}- \tilde{\boldsymbol \beta})^T{\bf X}^* \right\}^2=E\left[ (\alpha -\tilde{\alpha}) +  ({\boldsymbol \beta}- \tilde{\boldsymbol \beta})^T E({\bf X}) +  ({\boldsymbol \beta}- \tilde{\boldsymbol \beta})^T\{{\bf X}^*-E({\bf X})\} \right]^2\\
= E\left\{ (\alpha -\tilde{\alpha}) +  ({\boldsymbol \beta}- \tilde{\boldsymbol \beta})^T E({\bf X}) \right\}^2 + E\{ ({\boldsymbol \beta}- \tilde{\boldsymbol \beta})^T {\bf M} ({\boldsymbol \beta}- \tilde{\boldsymbol \beta})  \},
\end{multline}
For the first term in \eqref{eq:term1_} notice that 
$\alpha = E(Y) -\sum_{j=1}^p \beta_j E(X_j)$ and $\tilde{\alpha} = \bar{Y} - \sum_{j=1}^p \tilde{\beta}_j {\bar X}_j$; therefore,
\[
E\left\{ (\alpha -\tilde{\alpha}) +  ({\boldsymbol \beta}- \tilde{\boldsymbol \beta})^T E({\bf X}) \right\}^2 = E\Big[  \bar{Y} - \sum_{j=1}^p \tilde{\beta}_j\{ {\bar X}_j - E(X_j)\} -E(Y)\Big]^2,
\] 
and \eqref{eq:excess_risk_} is established. \qed

\subsection*{\bf Proof of Proposition \ref{prop1}.}

{\bf Proposition 1.}	Let  $\tilde{\alpha},\tilde{\boldsymbol \beta}$ satisfies 
	\begin{equation*}
	\tilde{\alpha} = \bar{Y} - \sum_{j=1}^p \tilde{\beta}_j {\bar X}_j,
	\end{equation*}
	and assume that $\tilde{\boldsymbol \beta}$ satisfies $\sqrt{n} (\tilde{\boldsymbol \beta}- {\boldsymbol \beta}) \tendd N(0 , {\boldsymbol \Sigma})$ then,
	\begin{equation*}
	{R}(\tilde{\alpha},\tilde{\boldsymbol \beta}) =  E\Big[  Y - \sum_{j=1}^p {\beta}_j\{ { X}_j - E(X_j)\} -E(Y)\Big]^2
	+ {\rm Trace} ({\bf M} {\boldsymbol \Sigma}),
	\end{equation*}
	where ${\bf M}$ is the covariance matrix of ${\bf X}$.

{\bf Proof:}
Since for any numbers $c,d,B$ we have that $\min(c+d,B) \le \min(c,B) + \min(d,B)$ and $\min(c+d,B) \ge \min(c,B/2) + \min(d,B/2)$, then
\begin{multline} \label{eq:tmm}
{R}(\tilde{\alpha}, \tilde{\boldsymbol \beta}) = \lim_{B \to \infty} \lim_{n \to \infty}  E\min\{ n \tilde{L}(\tilde{\alpha},\tilde{\boldsymbol \beta}), B\} =
\lim_{B \to \infty} \lim_{n \to \infty} E\min\left( n\Big[  \bar{Y} - \sum_{j=1}^p \tilde{\beta}_j\{ {\bar X}_j - E(X_j)\} -E(Y)\Big]^2,B\right) \\+  \lim_{B \to \infty} \lim_{n \to \infty}  E  \min\left(  n\left\{ ({\boldsymbol \beta}- \tilde{\boldsymbol \beta})^T {\bf M} ({\boldsymbol \beta}- \tilde{\boldsymbol \beta}) \right \},B\right).
\end{multline} 
We start with the first summand in \eqref{eq:tmm}. Since $\sqrt{n}(\tilde{\boldsymbol \beta}-{\boldsymbol \beta})$ is asymptotically normal then
\[
\bar{Y} - \sum_{j=1}^p \tilde{\beta}_j\{ \bar{ X}_j - E(X_j)\} -E(Y) - \Big[ \bar{Y} - \sum_{j=1}^p {\beta}_j\{ \bar{ X}_j - E(X_j)\} -E(Y) \Big] = 
\sum_{j=1}^p ({\beta}_j- \tilde{\beta}_j)\{\bar { X}_j - E(X_j)\}= o_p(1/\sqrt{n}). 
\]
Therefore, by the continuous mapping theorem and since the random variables are bounded, then for any $B$,
\begin{multline}\label{eq:tmm1}
\lim_{n \to \infty}  E \min\left( n\Big[  \bar{Y} - \sum_{j=1}^p \tilde{\beta}_j\{ {\bar X}_j - E(X_j)\} -E(Y)\Big]^2,B\right) \\
= \lim_{n \to \infty}  E \min\left( n\Big[  \bar{Y} - \sum_{j=1}^p {\beta}_j\{ {\bar X}_j - E(X_j)\} -E(Y)\Big]^2,B\right)\\
=E \min\left( Z_1^2,B\right),
\end{multline}
where $Z_1$ is normal with mean zero and variance $E\Big[  {Y} - \sum_{j=1}^p {\beta}_j\{ {X}_j - E(X_j)\} -E(Y)\Big]^2$.

For the second summand in \eqref{eq:tmm}, the asymptotic normality $\sqrt{n}(\tilde{\boldsymbol \beta}-{\boldsymbol \beta}) \tendd N(0, {\boldsymbol \Sigma})$ implies that  for any $B$,
\begin{equation}\label{eq:tmm2}
\lim_{n \to \infty} E\min\left( n \left\{ ({\boldsymbol \beta}- \tilde{\boldsymbol \beta})^T {\bf M} ({\boldsymbol \beta}- \tilde{\boldsymbol \beta}) \right \},B\right)=E \min\left\{ {\bf Z}_2^T {\bf Z}_2 ,B \right\},
\end{equation}
where ${\bf Z}_2$ is a normal vector with mean zero and variance matrix ${\bf M}^{1/2} {\boldsymbol \Sigma}{\bf M}^{1/2}$.
Taking limits as $B \to \infty$ in \eqref{eq:tmm1} and \eqref{eq:tmm2} completes the proof of the proposition since ${\rm Trace}({\bf M}^{1/2} {\boldsymbol \Sigma}{\bf M}^{1/2})={\rm Trace}({\bf M} {\boldsymbol \Sigma})$.
\qed

\subsection*{\bf Proof of Theorem \ref{thm:p_dim}}

 We will not prove Theorem \ref{thm:one_dim} since it is a special case of Theorem \ref{thm:p_dim}, which is stated and proved now.
 
{\bf Theorem 2.} \begin{enumerate}
	\item[(i)] Suppose that the vector ${\bf X}_{\bullet}=\left\{\frac{X_{1\bullet}}{E(X_{1 \bullet}^2)},\ldots,\frac{X_{p\bullet}}{E(X_{p \bullet}^2)}\right\}$ is well defined (i.e., the projections in \eqref{eq:proj} exist and $E(X_{j \bullet})^2$ is positive and finite for $j=1,\ldots,p$) and that  $\delta {\bf X}_{\bullet}$ has finite second moments,  then,
	\[
	\sqrt{n}(\hat{\boldsymbol \beta}_{LSE}-{\boldsymbol \beta}) \tendd N(0,{\boldsymbol \Sigma}_{LSE}),
	\]
	where 
	${\boldsymbol \Sigma}_{LSE}= Cov(\delta {\bf X}_{\bullet})$.
	\item[(ii)]  Suppose further that for each $j=1,\ldots,p$, the vector $(W_j,U_1,U_2,\ldots,U_{p+1})$ has finite second moments and that the matrix $E({\bf U}_j {\bf U}_j^T)$ is invertible for ${\bf U}_j= \left( 1, U_1, U_2,\ldots,U_{p+1} \right)^T$;
	then,
	\[
	\sqrt{n}(\hat{\boldsymbol \beta}_{TI}-{\boldsymbol \beta}) \tendd N(0,{\boldsymbol \Sigma}_{TI}),
	\]
	where ${\boldsymbol \Sigma}_{TI}= Cov(\tilde{\boldsymbol \delta})$ for $\tilde{\boldsymbol \delta}=(\tilde{\delta}_1,\ldots,\tilde{\delta}_p)$. Furthermore, we have that ${\boldsymbol \Sigma}_{LSE} =  {\boldsymbol \Sigma}_{TI}+{\boldsymbol \Sigma}_{\rm diff}$ where 
	${\boldsymbol \Sigma}_{\rm diff}=Cov({\bf V})$ for ${\bf V}=\delta {\bf X}_{\bullet}-\tilde{\boldsymbol \delta}$.
	\item[(iii)] Suppose further that $\lim \frac{n}{n+m} = \nu$, then, 
	\[
	\sqrt{n}(\hat{\boldsymbol \beta}_{PI}-{\boldsymbol \beta}) \tendd N(0,{\boldsymbol \Sigma_{PI}}),
	\]
	where ${\boldsymbol \Sigma}_{PI}= {\boldsymbol \Sigma}_{TI}+\nu{\boldsymbol \Sigma}_{\rm diff}$.
\end{enumerate}
Therefore, if ${\boldsymbol \Sigma}_{\rm diff}$ is not the zero matrix, then ${\boldsymbol \Sigma}_{LSE}-{\boldsymbol \Sigma}_{TI}$ is positive definite, and if further $\nu<1$, then ${\boldsymbol \Sigma}_{LSE}-{\boldsymbol \Sigma}_{PI}$ is also positive definite.

{\bf Proof:}
\noindent \underline{\em Part (i)}  
 We have that \citep{Buja} $\{\hat{\boldsymbol \beta}_{LSE}\}_j = \beta_j + \frac{\langle {\boldsymbol \delta},{\bf X}_{j \bar{\bullet}} \rangle}{|| {\bf X}_{j \bar{\bullet}} ||^2}$. Furthermore,
\[
\{\hat{\boldsymbol \beta}_{LSE}\}_j = \beta_j +\frac{\langle {\boldsymbol \delta},{\bf X}_{j \bar{\bullet}} \rangle}{|| {\bf X}_{j \bar{\bullet}} ||^2} = \beta_j + \frac{\frac{1}{n} \sum_{i=1}^n { X}_{j \bar{\bullet}}^{(i)} \delta^{(i)} }{\frac{1}{n}\sum_{i=1}^n \left\{ { X}_{j \bar{\bullet}}^{(i)}\right\}^2 }=\beta_j+\frac{\frac{1}{n} \sum_{i=1}^n { X}_{j {\bullet}}^{(i)} \delta^{(i)} }{E(X_{j\bullet}^2)}+o_p(1/\sqrt{n}),
\]
since 
\[
\frac{\sqrt{n}}{n} \sum_{i=1}^n\delta^{(i)} (  { X}_{j \bar{\bullet}} ^{(i)} - { X}_{j\bullet}^{(i)}) =\frac{\sqrt{n}}{n} \sum_{i=1}^n\delta^{(i)} \left\{ \vec{\bf X}_{-j}^{(i)}\right\}^T\left(  \hat{\boldsymbol \beta}_{-j\bar{\bullet}}-{\boldsymbol \beta}_{-j{\bullet}}\right)  \tendp 0,\text{ and } \frac{1}{n} || {\bf X}_{j \bar{\bullet}} ||^2 \tendp E(X_{j\bullet}^2). 
\]
The results now follows from the multivariate CLT and Slutsky's theorem.\\
\noindent \underline{\em Part (ii)}
 For each $j$,
\begin{multline} \label{eq:TI1}
\{\hat{\boldsymbol \beta}_{TI}\}_j = \bar{W}_j -\hat{a} \bar{U}_1 -\sum_{j' =1}^p\hat{b}_{j'} \bar{U}_{j'+1}  = \bar{W}_j -{a} \bar{U}_1 -\sum_{j'=1}^p {b}_{j'} \bar{U}_{j'+1} +o_p(1/\sqrt{n})\\
= \beta_j+ \frac{1}{n}\sum_{i=1}^n \tilde{\delta}_j^{(i)} +o_p(1/\sqrt{n}).
\end{multline}
Again, the  multivariate CLT and Slutsky's theorem imply the result.

To prove the furthermore part, notice that for each $j$,
\begin{multline} \label{eq:diff_j}
V_j =\frac{\delta X_{j\bullet}}{E(X_{j\bullet}^2)}-\tilde{\delta}_j = \frac{\left(Y-\alpha - \sum_{j'=1}^p \beta_{j'} X_{j'} \right) X_{j\bullet}}{E(X_{j\bullet}^2)}-\frac{YX_{j\bullet}}{E(X_{j\bullet}^2)} + \beta_j - b_j +a\frac{X_{j\bullet}}{E(X_{j\bullet}^2)} + \sum_{j' =1}^pb_{j'} \frac{X_{j'}X_{j\bullet}}{E(X_{j\bullet}^2)}\\
=\beta_j - {b_j} + (a-\alpha)\frac{ X_{j\bullet}}{E(X_{j\bullet}^2)}  + \sum_{j'=1}^p(b_{j'}-\beta_{j'}) \frac{X_{j'}X_{j\bullet}}{E(X_{j\bullet}^2)}= (a-\alpha)U_1  + \sum_{j'=1}^p(b_{j'}-\beta_{j'}) U_{j'+1} 
\end{multline}
is a linear function of $U_1,\ldots,U_{p+1}$ and hence is orthogonal to $\tilde{\delta}_j$, i.e., $EV_j \tilde{\delta}_j=0$ for all $j$. 
Therefore,
\[
{\boldsymbol \Sigma}_{LSE}=Cov( {\bf X}_{\bullet} \delta)=Cov\left\{ \tilde{\boldsymbol \delta} + ( {\bf X}_{\bullet} \delta - \tilde{\boldsymbol \delta}) \right\} =Cov( \tilde{\boldsymbol \delta} + {\bf V}) = Cov( \tilde{\boldsymbol \delta}) + Cov({\bf V})= {\boldsymbol \Sigma}_{TI}+{\boldsymbol \Sigma}_{\rm diff}.
\]
\noindent \underline{\em Part (iii)}
We show that 
\begin{equation}\label{eq:diff_npm}
\{\hat{\boldsymbol \beta}_{PI}\}_j - \{\hat{\boldsymbol \beta}_{TI}\}_j = \frac{1}{n+m} \sum_{i=1}^{n+m} V_j^{(i)} +o_p(1/\sqrt{n}). 
\end{equation}
 Hence, $\sqrt{n}(\hat{\boldsymbol \beta}_{PI} - \hat{\boldsymbol \beta}_{TI}) \tendd N(0, \nu {\boldsymbol \Sigma}_{\rm diff})$. Since ${\bf V}$ is orthogonal to $\tilde{\boldsymbol \delta}$, then
\[
\sqrt{n}(\hat{\boldsymbol \beta}_{PI} - {\boldsymbol \beta}) = \sqrt{n}(\hat{\boldsymbol \beta}_{TI} - {\boldsymbol \beta})+\sqrt{n}(\hat{\boldsymbol \beta}_{PI} - \hat{\boldsymbol \beta}_{TI}) \tendd N(0,{\boldsymbol \Sigma}_{TI} + \nu {\boldsymbol \Sigma}_{\rm diff}). 
\] 
 
We now prove \eqref{eq:diff_npm} by a somewhat lengthy calculation.
 We have that 
\begin{equation} \label{eq:PI1}
\{\hat{\boldsymbol \beta}_{PI}\}_j = \bar{\check W}_j -\hat{a} \bar{\check{U}}_1 -\sum_{j' =1}^p\hat{b}_{j'} \bar{\check{U}}_{j'+1}  = \bar{\hat W}_j -{a} \bar{\check{U}}_1 -\sum_{j' =1}^p {b}_{j'} \bar{\check{U}}_{j'+1} +o_p(1/\sqrt{n}).
\end{equation}
Therefore, \eqref{eq:PI1} and \eqref{eq:TI1} yield
\[
\{\hat{\boldsymbol \beta}_{PI}\}_j-\{\hat{\boldsymbol \beta}_{TI}\}_j=\bar{\check W}_j - \bar{W}_j -{a}( \bar{\check{U}}_1 - \bar{U}_1) -\sum_{j'=1}^p {b}_{j'} (\bar{\check{U}}_{j'+1}-\bar{U}_{j'+1}) +o_p(1/\sqrt{n}).
\]
We have that 
\[
\bar{\check W}_j - \bar{W}_j = {\alpha}( \bar{\check{U}}_1 - \bar{U}_1) +\sum_{j'=1}^p {\beta}_{j'} (\bar{\check{U}}_{j'+1}-\bar{U}_{j'+1}) +o_p(1/\sqrt{n}),
\]
and hence,
\begin{equation}\label{eq:term}
\{\hat{\boldsymbol \beta}_{PI}\}_j-\{\hat{\boldsymbol \beta}_{TI}\}_j =
(\alpha - a)( \bar{\check{U}}_1 - \bar{U}_1) +\sum_{j'=1}^p ({\beta}_{j'}-b_{j'}) (\bar{\check{U}}_{j'+1}-\bar{U}_{j'+1}) +o_p(1/\sqrt{n}).
\end{equation} 
We now consider the summands in \eqref{eq:term}. We start with the first summand
\begin{multline*}
 \bar{\check{U}}_1 - \bar{U}=\frac{\frac{1}{n}\sum_{i=1}^n X_{j\check{\bullet}}^{(i)}}{\check{E}(X_{j\check{\bullet}}^2)} -\frac{\frac{1}{n}\sum_{i=1}^n X_{j{\bullet}}^{(i)}}{{E}(X_{j{\bullet}})^2} = \frac{1}{n}\sum_{i=1}^n X_{j{\bullet}}^{(i)}\left( \frac{1}{\check{E}(X_{j\check{\bullet}}^2)} -\frac{1}{{E}(X_{j{\bullet}})^2}
\right) + \frac{({\boldsymbol \beta}_{-j {\bullet}}-{\boldsymbol \beta}_{-j \check{\bullet}})^T\frac{1}{n}\sum_{i=1}^n \vec{\bf X}_{-j}^{(i)}}{\check{E}(X_{j\check{\bullet}}^2)}\\
=  \frac{({\boldsymbol \beta}_{-j {\bullet}}-{\boldsymbol \beta}_{-j \check{\bullet}})^T\frac{1}{n}\sum_{i=1}^n \vec{\bf X}_{-j}^{(i)}}{{E}(X_{j{\bullet}})^2}+o_p(1/\sqrt{n})\\
=  \frac{({\boldsymbol \beta}_{-j {\bullet}}-{\boldsymbol \beta}_{-j \check{\bullet}})^T\frac{1}{n+m}\sum_{i=1}^{n+m} \vec{\bf X}_{-j}^{(i)}}{{E}(X_{j{\bullet}})^2}+\frac{({\boldsymbol \beta}_{-j {\bullet}}-{\boldsymbol \beta}_{-j \check{\bullet}})^T\left(\frac{1}{n}\sum_{i=1}^n \vec{\bf X}_{-j}^{(i)}-\frac{1}{n+m}\sum_{i=1}^{n+m} \vec{\bf X}_{-j}^{(i)}\right)}{{E}(X_{j{\bullet}})^2}+o_p(1/\sqrt{n})\\
=  \frac{({\boldsymbol \beta}_{-j {\bullet}}-{\boldsymbol \beta}_{-j \check{\bullet}})^T\frac{1}{n+m}\sum_{i=1}^{n+m} \vec{\bf X}_{-j}^{(i)}}{{E}(X_{j{\bullet}})^2}+o_p(1/\sqrt{n}).
\end{multline*}
Because
\[
0=\check{E}(X_{j \check{\bullet}} ) =\check{E}(X_j -  {\boldsymbol \beta}_{-j \check{\bullet}}^T \vec{\bf X}_{-j}) \Longrightarrow \check{E} (X_j)=\check{E}( {\boldsymbol \beta}_{-j \check{\bullet}}^T \vec{\bf X}_{-j}),
\] 
we have,
\[
\check{E} (X_{j {\bullet}} ) =  \check{E} (X_j -  {\boldsymbol \beta}_{-j \check{\bullet}}^T \vec{\bf X}_{-j})=\check{E}( {\boldsymbol \beta}_{-j \check{\bullet}}^T \vec{\bf X}_{-j})-\check{E}( {\boldsymbol \beta}_{-j {\bullet}}^T \vec{\bf X}_{-j})=
{({\boldsymbol \beta}_{-j \check{\bullet}} -{\boldsymbol \beta}_{-j {\bullet}})^T\frac{1}{n+m}\sum_{i=1}^{n+m} \vec{\bf X}_{-j}^{(i)}},
\]
and therefore,
\begin{equation}\label{eq:term1}
\bar{\check{U}}_1 - \bar{U}_1 = \frac{-\frac{1}{n+m} \sum_{i=1}^{n+m} X_{j {\bullet}}^{(i)}}{{E}(X_{j{\bullet}})^2}+o_p(1/\sqrt{n}).
\end{equation}

We now consider the $j'$-th term in \eqref{eq:term} for $j' \ne j$; then, $E(X_{j'}X_{j{\bullet}}) = 0$. We have
\begin{multline*}
 \bar{\check{U}}_{j'+1}-\bar{U}_{j'+1}=\frac{\frac{1}{n}\sum_{i=1}^n X^{(i)}_{j'}X^{(i)}_{j\check{\bullet}}}{\check{E}(X_{j\check{\bullet}}^2)}-\frac{\frac{1}{n}\sum_{i=1}^n X^{(i)}_{j'}X^{(i)}_{j{\bullet}}}{{E}(X_{j{\bullet}})^2}\\
=\frac{1}{n}\sum_{i=1}^n X^{(i)}_{j'}X^{(i)}_{j{\bullet}}\left\{\frac{1}{\check{E}(X_{j\check{\bullet}}^2)}-\frac{1}{{E}(X_{j{\bullet}})^2} \right\}+ \frac{({\boldsymbol \beta}_{-j {\bullet}}-{\boldsymbol \beta}_{-j \check{\bullet}})^T\frac{1}{n}\sum_{i=1}^n X_{j'}^{(i)} \vec{\bf X}_{-j}^{(i)}}{\check{E}(X_{j\check{\bullet}}^2)}\\
=  \frac{({\boldsymbol \beta}_{-j {\bullet}}-{\boldsymbol \beta}_{-j \check{\bullet}})^T\frac{1}{n+m}\sum_{i=1}^{n+m} X_{j'}^{(i)} \vec{\bf X}_{-j}^{(i)}}{{E}(X_{j{\bullet}})^2}+o_p(1/\sqrt{n}).
\end{multline*}
Because
\[
0=\check{E}(X_{j \check{\bullet}} X_{j'}) =\check{E}\{ (X_j -  {\boldsymbol \beta}_{-j \check{\bullet}}^T \vec{\bf X}_{-j})X_{j'}\} \Longrightarrow \check{E} (X_j X_{j'})=\check{E}( {\boldsymbol \beta}_{-j \check{\bullet}}^T \vec{\bf X}_{-j}X_{j'}),
\] 
we have,
\[
\check{E} (X_{j {\bullet}} X_{j'} ) =  \check{E}\{ (X_j -  {\boldsymbol \beta}_{-j \check{\bullet}}^T \vec{\bf X}_{-j})X_{j'}\}=\check{E}( {\boldsymbol \beta}_{-j \check{\bullet}}^T \vec{\bf X}_{-j}X_{j'})-\check{E}( {\boldsymbol \beta}_{-j {\bullet}}^T \vec{\bf X}_{-j}X_{j'})=
{({\boldsymbol \beta}_{-j \check{\bullet}} -{\boldsymbol \beta}_{-j {\bullet}})^T\frac{1}{n+m}\sum_{i=1}^{n+m} X_{j'}^{(i)} \vec{\bf X}_{-j}^{(i)}},
\]
and therefore,
\begin{equation}\label{eq:term2}
\bar{\check{U}}_{j'+1}-\bar{U}_{j'+1} = \frac{-\frac{1}{n+m} \sum_{i=1}^{n+m} X_{j {\bullet}}^{(i)} X_{j'}^{(i)}}{{E}(X_{j{\bullet}})^2}+o_p(1/\sqrt{n}).
\end{equation}

For the $j$-th term in \eqref{eq:term}, we have that $E(X_{j}X_{j{\bullet}}) = E(X_{j{\bullet}})^2$; hence,
\begin{multline*}
\bar{\check{U}}_{j+1}-\bar{U}_{j+1}= \frac{\frac{1}{n}\sum_{i=1}^n X^{(i)}_{j}X^{(i)}_{j\check{\bullet}}}{\check{E}(X_{j\check{\bullet}}^2)}-\frac{\frac{1}{n}\sum_{i=1}^n X^{(i)}_{j}X^{(i)}_{j{\bullet}}}{{E}(X_{j{\bullet}})^2}\\
=\left\{ \frac{1}{n}\sum_{i=1}^n X^{(i)}_{j'}X^{(i)}_{j{\bullet}}-E(X_{j{\bullet}})^2\right\}\left\{\frac{1}{\check{E}(X_{j\check{\bullet}}^2)}-\frac{1}{{E}(X_{j{\bullet}})^2} \right\}
+ \frac{({\boldsymbol \beta}_{-j {\bullet}}-{\boldsymbol \beta}_{-j \check{\bullet}})^T\frac{1}{n}\sum_{i=1}^n X_{j}^{(i)} \vec{\bf X}_{-j}^{(i)}+E(X_{j{\bullet}})^2-\check{E}(X_{j\check{\bullet}}^2)}{\check{E}(X_{j\check{\bullet}}^2)}\\
=  \frac{({\boldsymbol \beta}_{-j {\bullet}}-{\boldsymbol \beta}_{-j \check{\bullet}})^T\frac{1}{n}\sum_{i=1}^n X_{j}^{(i)} \vec{\bf X}_{-j}^{(i)}+E(X_{j{\bullet}})^2-\check{E}(X_{j\check{\bullet}}^2)}{{E}(X_{j{\bullet}})^2}+o_p(1/\sqrt{n}).
\end{multline*}
Because 
\[
\check{E}(X_{j\check{\bullet}}^2)=\check{E}(X_{j \check{\bullet}} X_{j}) =\check{E}\{ (X_j -  {\boldsymbol \beta}_{-j \check{\bullet}}^T \vec{\bf X}_{-j})X_{j}\} \Longrightarrow \check{E} (X_j ^2)=\check{E}( {\boldsymbol \beta}_{-j \check{\bullet}}^T \vec{\bf X}_{-j}X_{j})+\check{E}(X_{j\check{\bullet}}^2),
\] 
we have,
\begin{multline*}
\check{E} (X_{j {\bullet}} X_{j} ) =  \check{E}\{ (X_j -  {\boldsymbol \beta}_{-j \check{\bullet}}^T \vec{\bf X}_{-j})X_{j}\}=\check{E}( {\boldsymbol \beta}_{-j \check{\bullet}}^T \vec{\bf X}_{-j}X_{j})+\check{E}(X_{j\check{\bullet}}^2)-\check{E}( {\boldsymbol \beta}_{-j {\bullet}}^T \vec{\bf X}_{-j}X_{j})\\
={({\boldsymbol \beta}_{-j \check{\bullet}} -{\boldsymbol \beta}_{-j {\bullet}})^T\frac{1}{n+m}\sum_{i=1}^{n+m} X_{j}^{(i)} \vec{\bf X}_{-j}^{(i)}}+\check{E}(X_{j\check{\bullet}}^2).
\end{multline*}
Thus, for the $j$-th term we obtain,
\begin{equation}\label{eq:term3}
 \bar{\check{U}}_{j+1}-\bar{U}_{j+1}=\frac{-\frac{1}{n+m} \sum_{i=1}^{n+m} X_{j {\bullet}}^{(i)} X_{j}^{(i)}}{{E}(X_{j{\bullet}})^2}+1+o_p(1/\sqrt{n}). 
\end{equation}

Tracking back to \eqref{eq:term} through \eqref{eq:term1},\eqref{eq:term2}, and \eqref{eq:term3}, we have
\[
\{\hat{\boldsymbol \beta}_{PI}\}_j-\{\hat{\boldsymbol \beta}_{TI}\}_j =(\beta_j- {b_j} )+ (a-\alpha)\frac{\frac{1}{n+m} \sum_{i=1}^{n+m} X_{j {\bullet}}^{(i)}}{E(X_{j\bullet}^2)}  + \sum_{j'=1}^p(b_{j'}-\beta_{j'})\frac{\frac{1}{n+m} \sum_{i=1}^{n+m} X_{j {\bullet}}^{(i)} X_{j'}^{(i)}}{E(X_{j\bullet}^2)}+o_p(1/\sqrt{n}).
\]
Hence, \eqref{eq:diff_j} implies that 
\begin{equation*}
\{\hat{\boldsymbol \beta}_{PI}\}_j-\{\hat{\boldsymbol \beta}_{TI}\}_j = \frac{1}{n+m} \sum_{i=1}^{n+m} V_j^{(i)} +o_p(1/\sqrt{n}), 
\end{equation*}
and \eqref{eq:diff_npm} is established.
\qed

\subsection*{Proof of Theorem \ref{thm:min}}
\noindent{\bf Theorem 3.} (Sufficiency) If $X \in \text{span}\{q_1(X),\ldots,q_K(X)\}$ and $X^2-1 \in \text{span}\{q_1(X),\ldots,q_K(X)\}$, then $\hat{\beta}_{K}$ is better than $\hat{\beta}_{LSE}$  for ``every'' $G$, i.e., $E_G\tilde{\delta}_K^2 \le E_G(\delta X)^2$ for every distribution $G$ that satisfies \eqref{eq:cond} and \eqref{eq:cond1}.\\
(Necessity) If $X \notin \text{span}\{q_1(X),\ldots,q_K(X)\}$ or $X^2-1 \notin \text{span}\{q_1(X),\ldots,q_K(X)\}$, then there exists distribution $G$ that satisfies \eqref{eq:cond} and \eqref{eq:cond1}, and where $\hat{\beta}_{LSE}$ is better than $\hat{\beta}_{K}$, i.e., $E_G(\delta X)^2 < E_G\tilde{\delta}_K^2$.

{\bf Proof:} (Sufficiency). The same proof outlined in the introduction applies. We bring it here for completeness. We have that
\[
\delta X - \tilde{\delta}_K =  \sum_{k=1}^K b_k q_k(X)- \alpha X - \beta (X^2-1). 
\]
If $X \in \text{span}\{q_1(X),\ldots,q_K(X)\}$ and $X^2-1 \in \text{span}\{q_1(X),\ldots,q_K(X)\}$, then $\delta X - \tilde{\delta}_K$ is orthogonal to $\tilde{\delta}_K$ and then
\[
E_G(\delta X)^2= E_G \left( \tilde{\delta}_K + \delta X - \tilde{\delta}_K \right)^2 = E_G \tilde{\delta}_K^2 + E_G \left(  \delta X - \tilde{\delta}_K \right)^2 \ge E_G \tilde{\delta}_K^2. 
\]  
(Necessity). Let $G$ be such that $X$ and $Y$ are linearly dependent, i.e., $Y=\alpha + \beta X$ for certain $\alpha$ and $\beta$. Then $\delta \equiv 0$. We have that 
\[
XY=\beta+\alpha X + \beta (X^2-1) 
\]
and therefore $\tilde{\delta}_K \equiv 0$ in Model \eqref{eq:int_model} for every $\alpha,\beta$ iff $X \in \text{span}\{q_1(X),\ldots,q_K(X)\}$ and $X^2-1 \in \text{span}\{q_1(X),\ldots,q_K(X)\}$.\qed

\section{Details on simulation results}\label{sec:sims}
\small

\begin{enumerate}
	\item{ A proposed rule of thumb when $p = 1$: As indicated by Table \ref{tab:prop} across a wide range of parameter settings, PI performs almost as well as LSE when $n \geq 100$, $m \geq 100$ and linear $E[Y|X]$, and better than the LSE for nonlinear $E[Y|X]$ (that is, $ERR < 1$). Furthermore, PI's underperformance is mild when $E[Y|X]$ is linear. Thus, assuming that the scenarios we studied are representative of reality, we recommend using the PI method when the following conditions are satisfied, as a rule of thumb:
		\begin{enumerate}
			\item{$n \geq 100$}
			\item{$m \geq 100$}
			\item{$E[Y|X]$ could be non-linear}
		\end{enumerate}
	}
	\item{If $p = 1$, $n > 100$, $m > 100$, and $E[Y|X]$ is one of the non-linear functions tested, then both PI and TI appear to outperform LSE; furthermore, the margin of outperformance widens as $n$ increases (See Figure \ref{fig:1} (a), and Figure \ref{fig:1} (b), respectively). In particular, PI does much better for $E[Y | {X}] = X^3$ and $e^X$, as compared to $\sqrt{X}$. Intuitively, this results from the larger non-linear moments of $X^3$ and $e^X$, causing ${\boldsymbol \Sigma}_{\rm diff}$ to be large (as defined in Theorem \ref{thm:p_dim}, part ii)}
	\item{Increasing $m$, the size of the unlabeled dataset, improves performance of PI. Indeed, comparing PI where $m = 2n$ in Graph X against TI ($m \geq 500n$) in Graph Y, one observes that greatly increasing the pool of unlabeled samples can reduce Excess Risk, in some cases on the order of $5-10$ per cent (Compare Figures \ref{fig:1} (a) and (b)). But, a word of caution: if $n$ is small, then LSE may still outperform both PI and TI, even when $m$ is large and $E[Y|X]$ is non-linear.}
	\begin{figure}[ht!] 
		\subfigure[PI ]{%
			\includegraphics[width=0.5\textwidth]{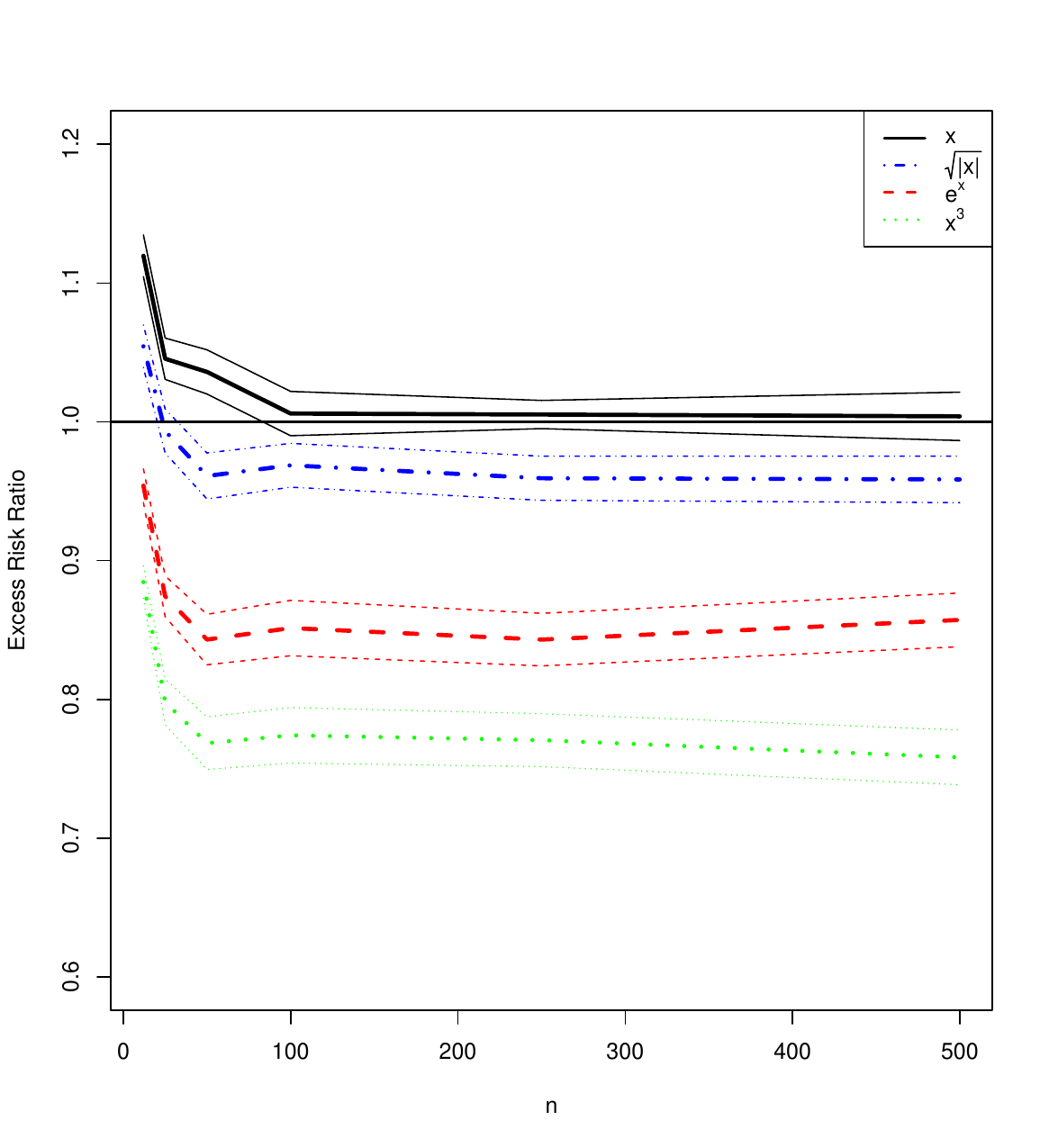}
		}
		\subfigure[TI]{%
			\includegraphics[width=0.5\textwidth]{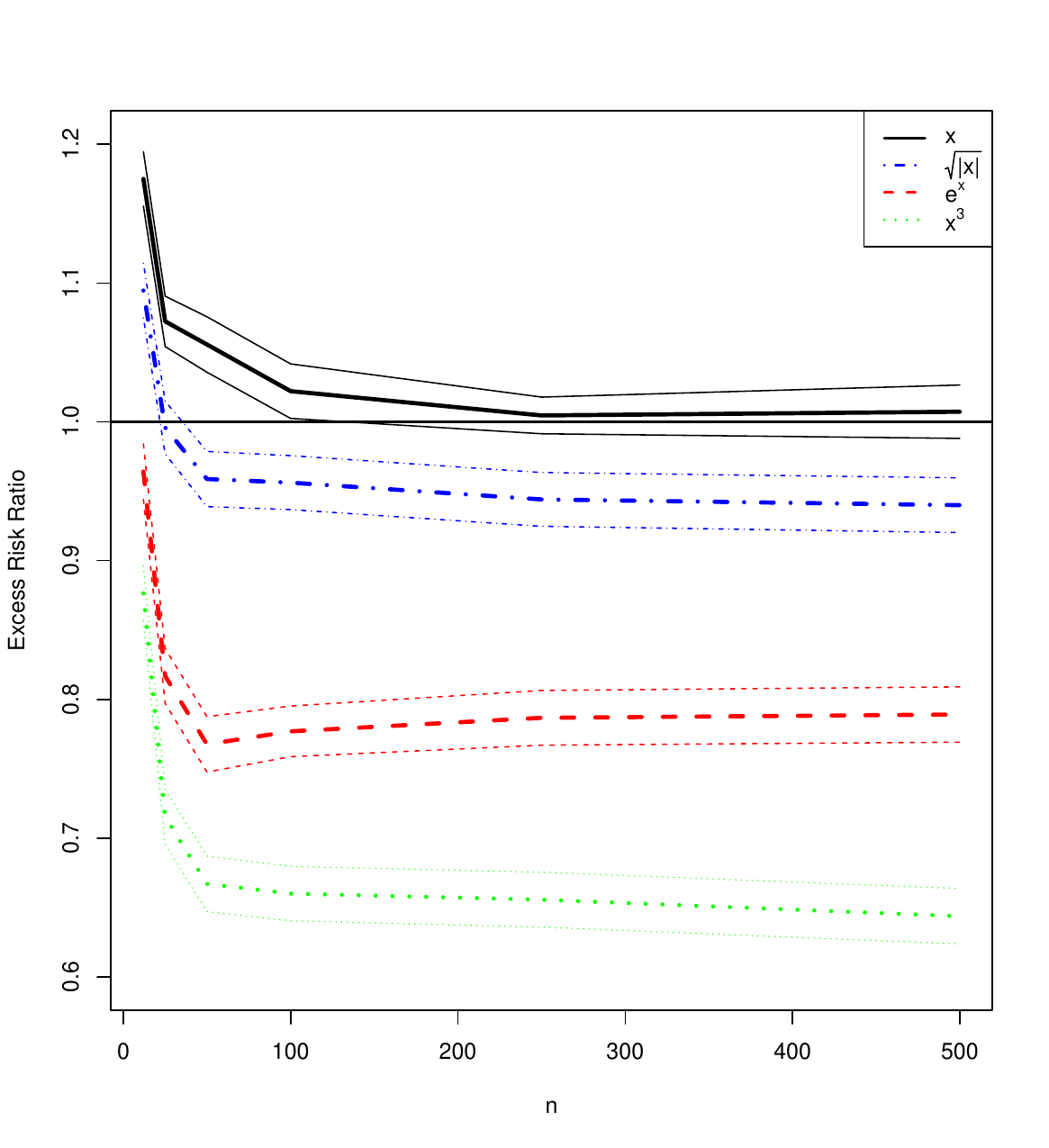}
		}
		\\    \caption{\footnotesize Estimates of the excess risk ratio for the TI and PI estimates, for several possible functions of $E[Y|X]$.   
			In PI estimates (a), $m = 2n$. In TI estimates (b), the PI method is used with $m \geq 500n$.
			In all cases, X is Gaussian, errors are Gaussian, and $p = 1$. Confidence intervals are $\pm 2$ standard errors. \label{fig:1}}
	\end{figure}     
	\item{If $n$ is small, or if $E[Y|X]$ is linear, then LSE performs better than TI and PI. However, $ERR \to 1$ as $n \to \infty$, holding other parameters fixed (See Figures \ref{fig:1} (a) and (b)). This result is consistent with the asymptotic agreement of PI, TI, and LSE when $E[Y|X]$ is linear (See Theorem \ref{thm:p_dim})}
	\item{Extra heteroskedastic noise appears to decrease the relative advantage of PI. See Figure \ref{fig:2} (a).}	
	\item{Somewhat similarly, increasing the number of parameters to $4$ (with $E[Y|X]$ being the sum over $j$ of the given function of $X_j$) also decreases the relative advantage of PI. See Figure \ref{fig:2} (b).}
	\item{If $X$ is not Gaussian, PI still does well, especially for large $n$. See Figure \ref{fig:2} (c) for an example where $X$ is exponential.
		(Note that standard errors may be understated, especially in the case $E[Y|X] = \exp(X)$, partly due to the large moments of $Y$.)}
	\begin{figure}[p!]        
		\subfigure[Heteroskedastic Errors]{%
			\includegraphics[width=0.5\textwidth]{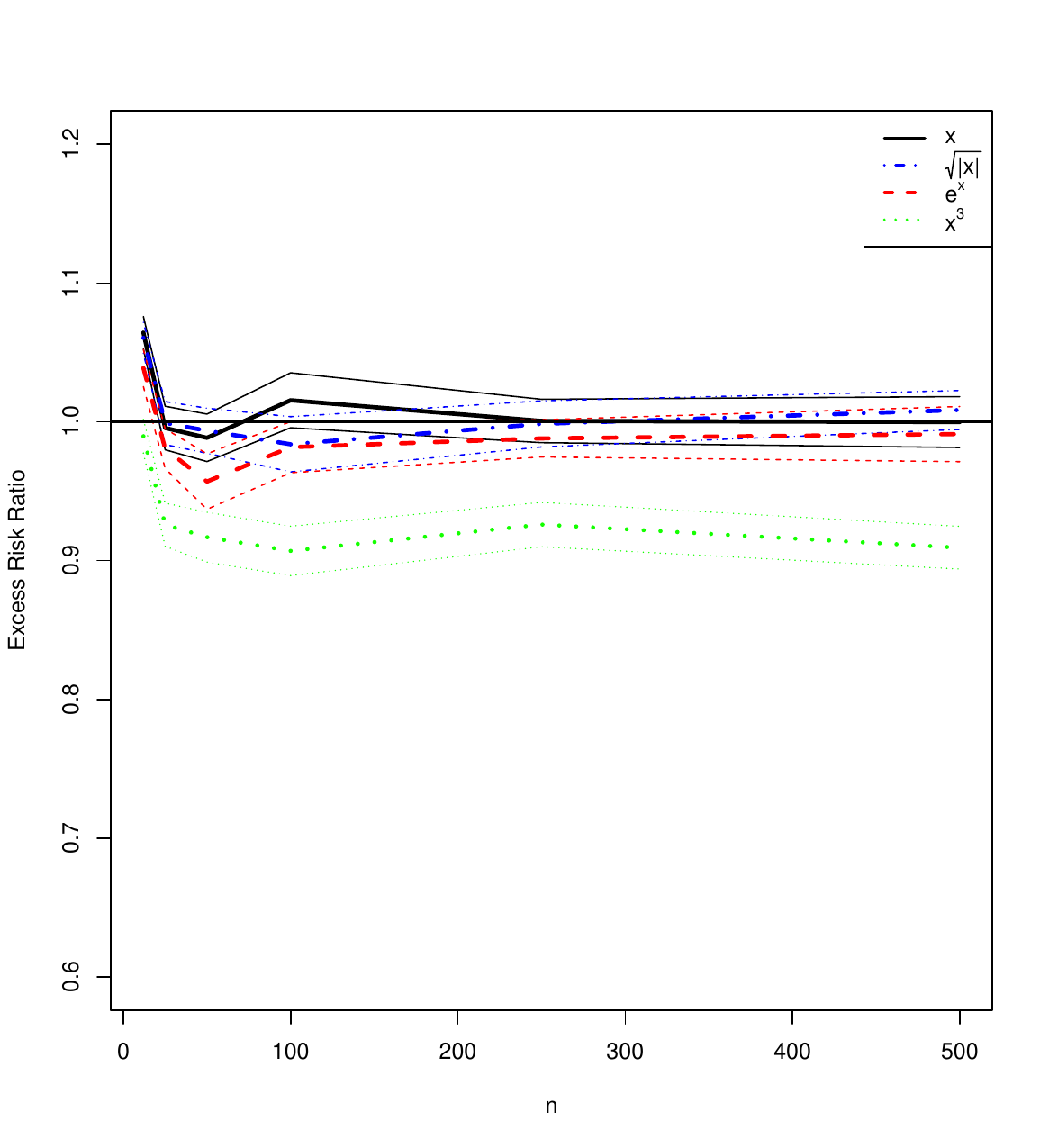}
		}
		\subfigure[$p = 4$]{%
			\includegraphics[width=0.5\textwidth]{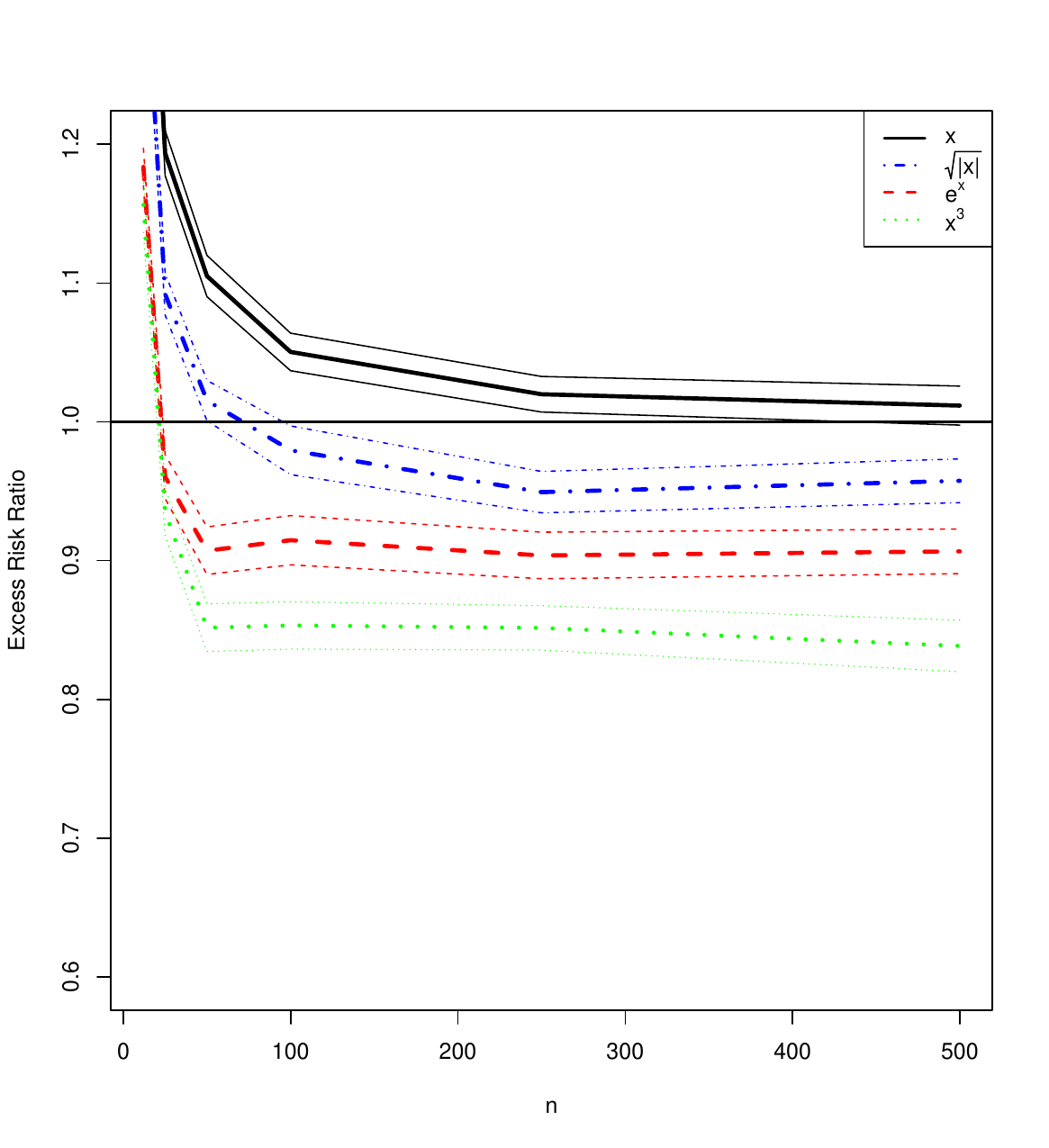}
		}
		\subfigure[$X$ Exponential ]{%
			\includegraphics[width=0.5\textwidth]{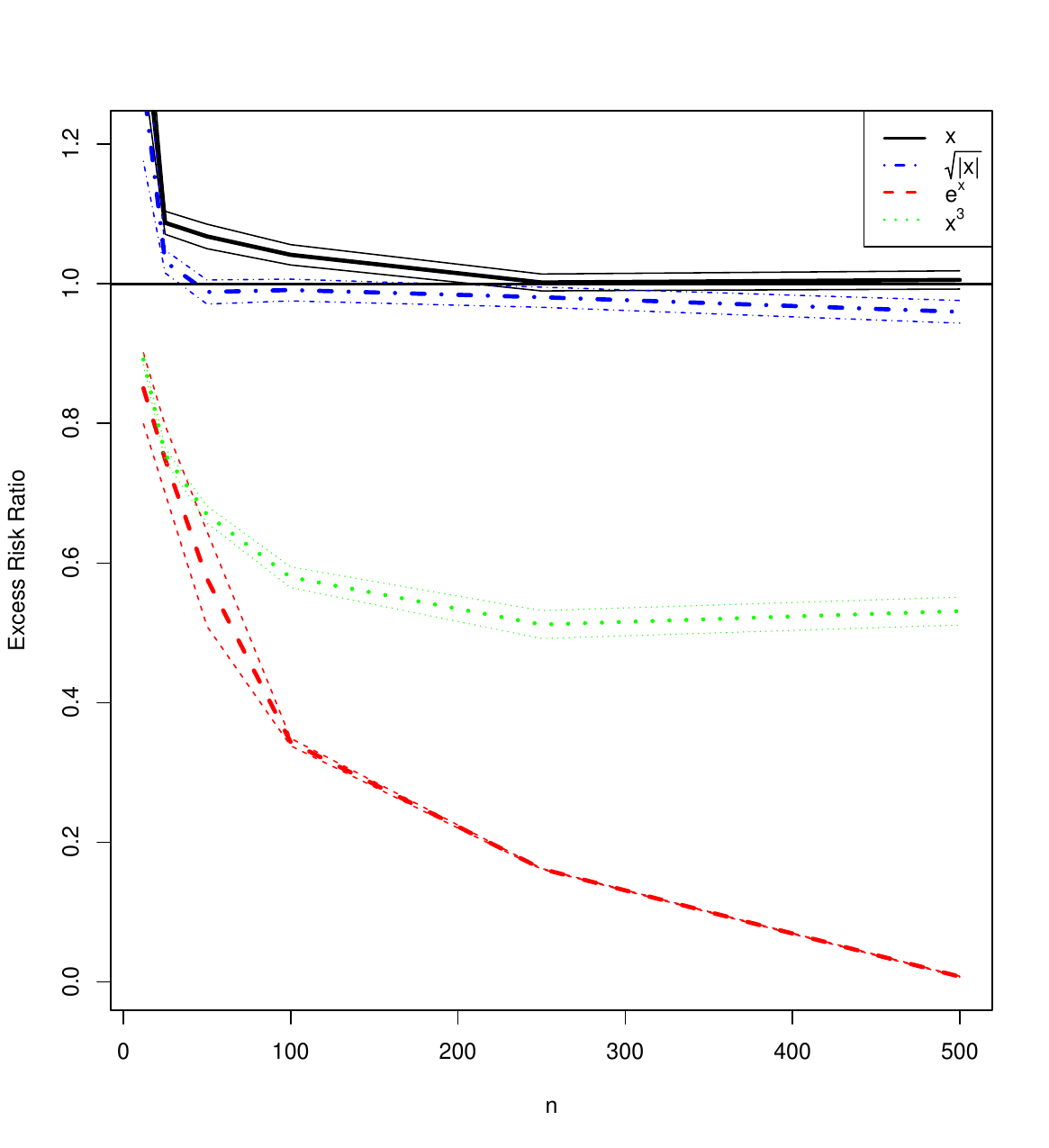}
		}
		\\    \caption{\footnotesize \linespread{1.3}\selectfont{} Estimates of the excess risk ratio of PI estimates, across several possible functions of $E[Y|X]$, with tweaks to the previous settings of Figure \ref{fig:1} (a): these were $m = 2n$, $X$ Gaussian, $p = 1$, and Guassian errors. The following is changed in each plot:
			Figure (a): Errors are heteroscedastic, distributed $N(0,e^{2\|x\|})$. Figure (b): $p = 4$, where $E[Y|X]$ is the sum of the given relationship in each $X$-variable. Figure (c): $X$ is distributed exponentially.
			Confidence intervals are $\pm 2$ standard errors.\label{fig:2}}
	\end{figure}		
\end{enumerate}

\section{List of the covariates for the bike sharing dataset}\label{sec:append}

The following covariates were used 
\begin{enumerate}
	\item Summer (binary).
	\item Fall (binary).
	\item Winter (binary).
	\item Hour (continuous).
	\item Holiday (binary).
	\item Sunday until Thursday (binary).
	\item Friday (binary).
	\item Weather: Mist + Cloudy, Mist + Broken clouds, Mist + Few clouds, Mist (binary).
	\item Temperature (continuous).
	\item Humidity (continuous).
\end{enumerate}

\section{Estimates of ${\bs \Sigma}_{PI}$}\label{sec:app_est}

In this section we study estimation of the asymptotic covariance matrix. We suggest three estimates for ${\boldsymbol \Sigma}_{PI}$ and study their performance by simulation. Estimation of the prediction error (excess risk) is also considered.

\subsection{Three estimates of ${\boldsymbol \Sigma}_{PI}$} \label{sec:est_var}

We proved that under certain conditions,
\[
\sqrt{n}(\hat{\boldsymbol \beta}_{LSE}-{\boldsymbol \beta}) \tendd N(0,{\boldsymbol \Sigma}_{LSE}),~\sqrt{n}(\hat{\boldsymbol \beta}_{PI}-{\boldsymbol \beta}) \tendd N(0,{\boldsymbol \Sigma}_{PI})\text{ and, }\sqrt{n}(\hat{\boldsymbol \beta}_{TI}-{\boldsymbol \beta}) \tendd N(0,{\boldsymbol \Sigma}_{TI}).
\]
In these expressions the (asymptotic) variability of the estimates are determined by the matrices ${\boldsymbol \Sigma}_{LSE}$, ${\boldsymbol \Sigma}_{PI}$ and ${\boldsymbol \Sigma}_{TI}$. These matrices are unknown and therefore need to be estimated. In this section we consider only estimation of ${\boldsymbol \Sigma}_{LSE}$ and ${\boldsymbol \Sigma}_{PI}$. 
Estimation of the asymptotic variance yields prediction error estimates or $L_2$ risk in a standard fashion. In section \ref{sec:pred_error} we discuss estimation of the excess risk.

For ${\boldsymbol \Sigma}_{LSE}$ we follow \citet{Buja} who show that the standard estimates are inconsistent when non-linearity is present and suggest two alternative estimates for the variance. The parametric sandwich estimator for the asymptotic variance of ${\hat {\boldsymbol \beta}}_{LSE}$ is 
\[
\hat{{\bs \Sigma}}_{LSE,PAR} = \left(\mathbb{X}^T {\mathbb X} \right)^{-1} \mathbb{X}^T \hat{D} {\mathbb X} \left(\mathbb{X}^T {\mathbb X} \right)^{-1},
\]
where $\mathbb{X}$ is the design matrix and $\hat{D}$ is a diagonal matrix with $\hat{D}_i=\hat{\delta}_i^2$ the standard residual estimator.
The second estimate is derived from a pairs bootstrap where a pair $({\bf X}^*,Y^*)$ is sampled $N_{BS}$ times from the empirical joint distribution of $({\bf X},Y)$, yielding a sample (under the empirical measure) of size $N_{BS}$ of ${\hat {\boldsymbol \beta}}_{LSE}$. The resulting estimate of the variance is denoted by $\hat{{\bs \Sigma}}_{LSE,BS}$.

For estimating the variance of ${\boldsymbol \Sigma}_{PI}$, we consider three estimates described as follows:
\begin{enumerate}
	\item {\bf Parametric:} 
	The asymptotic variances of ${\hat {\boldsymbol \beta}}_{LSE}$ and ${\hat {\boldsymbol \beta}}_{PI}$ are
	\begin{equation} \label{eq:AV_j}
	{{\bs \Sigma}}_{LSE} = Cov( \tilde{\boldsymbol \delta} ) +  Cov( {\bf  V} ), ~{{\bs \Sigma}}_{PI} = Cov( \tilde{\boldsymbol \delta}) +  \nu Cov( {\bf V} ).
	\end{equation}
	To estimate the variance of $Cov( \tilde{\boldsymbol \delta} )$ notice that for $j=1,\ldots,p$ and $i=1,\ldots,n$,
	\[
	\check{ \tilde{\delta}}_j^{(i)}= \check{W}^{(i)} -\{\hat{\boldsymbol \beta}_{PI}\}_j  - \hat{a} \check{U}_1^{(i)} -\sum_{j' =1}^p\hat{b}_{j'}\check{U}_{j'+1}^{(i)} =  \tilde{\delta}_j^{(i)} + o_p(1/\sqrt{n})  ,
	\]  
	and hence, $\overline{Cov}(  \check{ \tilde{\boldsymbol \delta}} ) \tendp Cov( \tilde{\boldsymbol \delta} )$, where $\overline{Cov}$ is the empirical covariance based on the labeled $n$ sample. Therefore, using \eqref{eq:AV_j} we obtain that a consistent estimator to ${{\bs \Sigma}}_{PI}$ is
	\[
	\left\{\hat{{\bs \Sigma}}_{PI,PAR}\right\}_{j,j'}= \overline{Cov} ( \check{ \tilde{\boldsymbol \delta}}_j, \check{ \tilde{\boldsymbol \delta}}_{j'}) + \nu  \left( \left\{\hat{{\bs \Sigma}}_{LSE,PAR}\right\}_{j,j'} -\overline{Cov} ( \check{ \tilde{\boldsymbol \delta}}_j, \check{ \tilde{\boldsymbol \delta}}_{j'})\right),  
	\]  
	for $j\ne j'$ and for $j=1,\ldots,p$,
	\[
	\left\{\hat{{\bs \Sigma}}_{PI,PAR}\right\}_{j,j}=\frac{1}{n} \sum_{i=1}^n \left\{ \check{ \tilde{\delta}}_j^{(i)} \right\}^2 + \nu \max \left(  \left\{\hat{{\bs \Sigma}}_{LSE,PAR}\right\}_{j,j} -\frac{1}{n} \sum_{i=1}^n \left\{ \check{ \tilde{\delta}}_j^{(i)} \right\}^2 ,0 \right).
	\]
	\item {\bf Bootstrap:} Let $\left\{ ({\bf X}^{(i)}_*,Y^{(i)}_*) \right\}_{i=1}^n$ be a sample from the empirical distribution of $({\bf X},Y)$ (based on the labeled observations) and $\{{\bf X}^{(i)}_*\}_{i=n+1}^{n+m}$ be a sample from the empirical distribution of ${\bf X}$ (based on the unlabeled observations). 
	Thus, ${\hat {\boldsymbol \beta}}_{PI}^* = {\hat {\boldsymbol \beta}}_{PI}^* \left( \left\{ ({\bf X}^{(i)}_*,Y^{(i)}_*) \right\}_{i=1}^n,\{{\bf X}^{(i)}_*\}_{i=n+1}^{n+m} \right)$ is a sample under the empirical measure of ${\hat {\boldsymbol \beta}}_{PI}$, which can yield an estimate of the variance of the estimator, denoted by $\hat{{\bs \Sigma}}_{PI,BS}$.    
	\item {\bf Variance bootstrap:}
	See Section \ref{sec:est_se}.
\end{enumerate}

In the next section we compare the different estimates through a simulation study. 

\subsection{Simulations for the estimates of the variance} \label{sec:est_var_sim} 

We now investigate the estimates of the variance under the toy model \eqref{eq:toymodel}. We study the performance of the estimates for $\alpha=1,1/2,1/4,1/8$, where small $\alpha$ corresponds to little non-linearity. We considered $n=250$ and $m=1500$ similar to the numbers of Section \ref{sec:inferring}. 
For each $\alpha$ we repeated the simulation 1000 times and also $N_{BS}=1000$. 

The simulations are summarized in Table \ref{tab:AV} and in Figure \ref{fig:AV}. We find that the bootstrap estimate (henceforth BS) is more variable than the parametric (henceforth PAR) and the variance bootstrap (henceforth VBS).
PAR and VBS are comparable but for small $\alpha$'s for estimating the difference, the histogram of VBS is narrower around the true value. When $\alpha=1/8$, the estimated difference under BS (PAR) is negative 32.9\% (34.7\%) of the simulations, while for VBS it is never negative. However, since VBS cannot be negative, the mean of the simulations is upwards biased. In short, we find that PAR and VBS outperform BS and VBS has the advantage of never being smaller than $\hat{{\bs \Sigma}}_{LSE,BS}$.

\begin{table} [ht!] 
	\caption{\footnotesize The mean (std) of the estimates of the asymptotic variance  ${\bs \Sigma}_{PI}$ and the difference ${\bs \Sigma}_{LSE}-{\bs \Sigma}_{PI}$ (here the matrices are one dimensional). }
	\begin{center}
		{
			\begin{tabular}{c|| c c c c |c cc c}
				{} &  \multicolumn{4}{c}{Variance} & \multicolumn{4}{c}{Difference}\\
				$\alpha$& Bootstrap &  Parametric & Variance BS & true  & Bootstrap &  Parametric & Variance BS & true \\
				\hline
				1 & 7.65 (3.71) & 6.86 (3.11) & 6.58 (2.40) & 7 & 2.63 (2.53) & 3.51 (1.53) & 3.70 (2.13) & 3.43 \\
				1/2 & 2.65 (0.97) & 2.42 (0.82) & 2.36 (0.70) & 2.5 & 0.64 (0.53) & 0.88 (0.45) & 0.93 (0.46) & 0.86 \\
				1/4 & 1.43 (0.40) & 1.33 (0.33) & 1.31 (0.31) & 1.37& 0.15 (0.18) & 0.24 (0.22) & 0.27 (0.16) & 0.21 \\
				1/8 & 1.11 (0.23) & 1.06 (0.20) & 1.03 (0.19) & 1.09 & 0.02 (0.07) & 0.06 (0.17) & 0.10 (0.06) & 0.05
		\end{tabular}}
		\label{tab:AV}
	\end{center}
\end{table}

\begin{figure}[ht!] 
	\subfigure[$\alpha=1/4$]{%
		\includegraphics[width=0.5\textwidth]{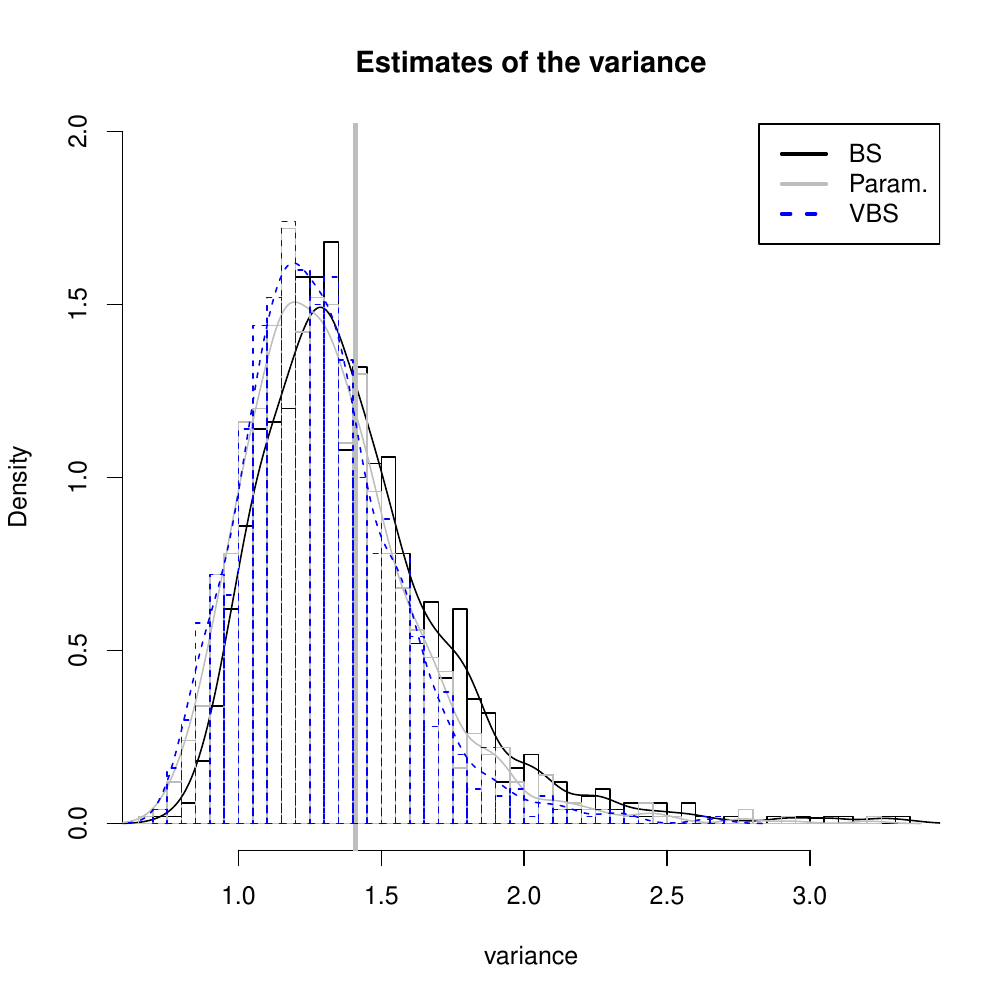}
		\includegraphics[width=0.5\textwidth]{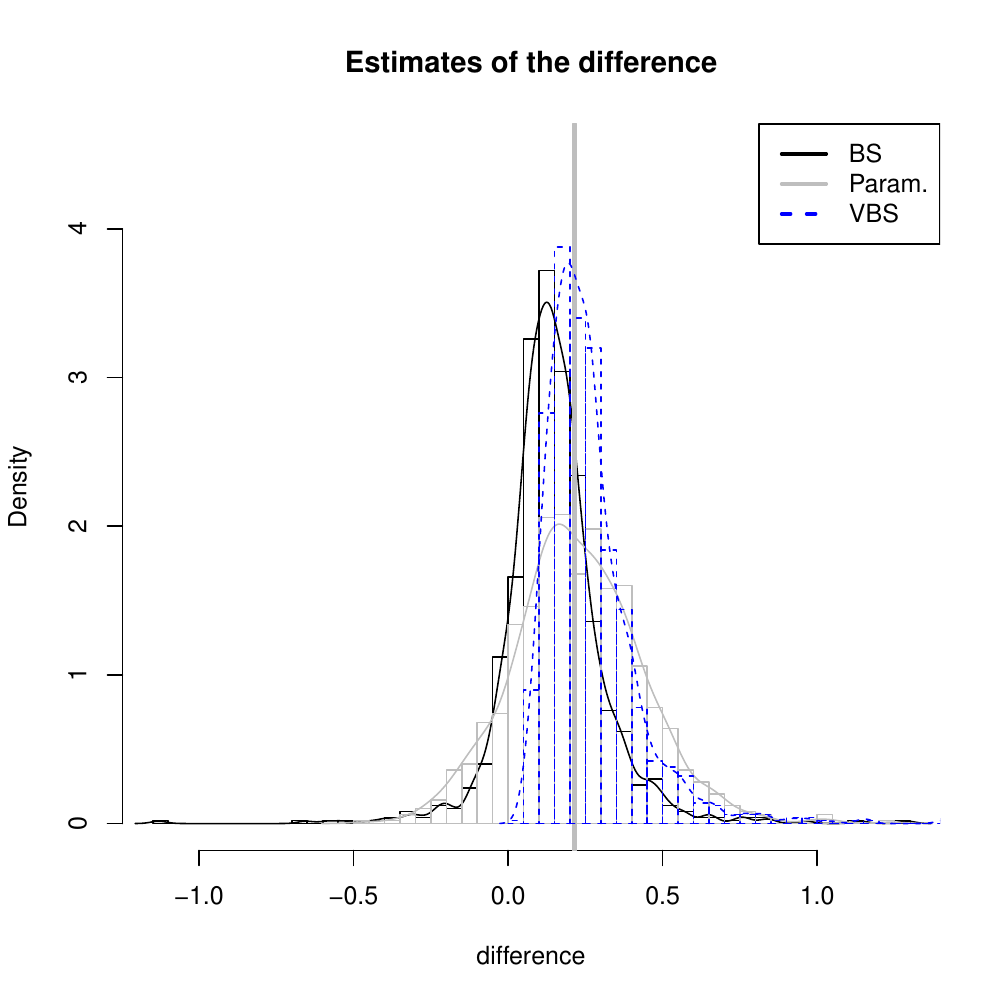}
	}\\
	\subfigure[$\alpha=1/8$]{%
		\includegraphics[width=0.5\textwidth]{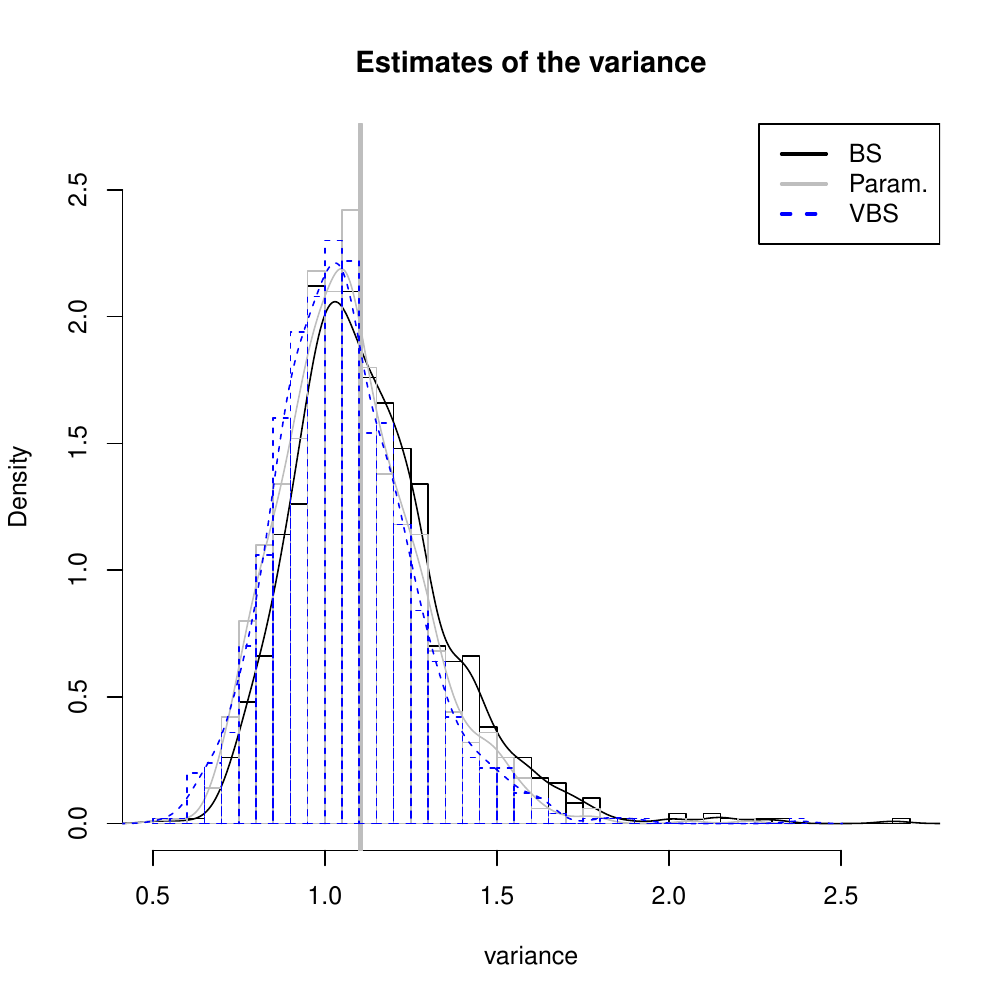}
		\includegraphics[width=0.5\textwidth]{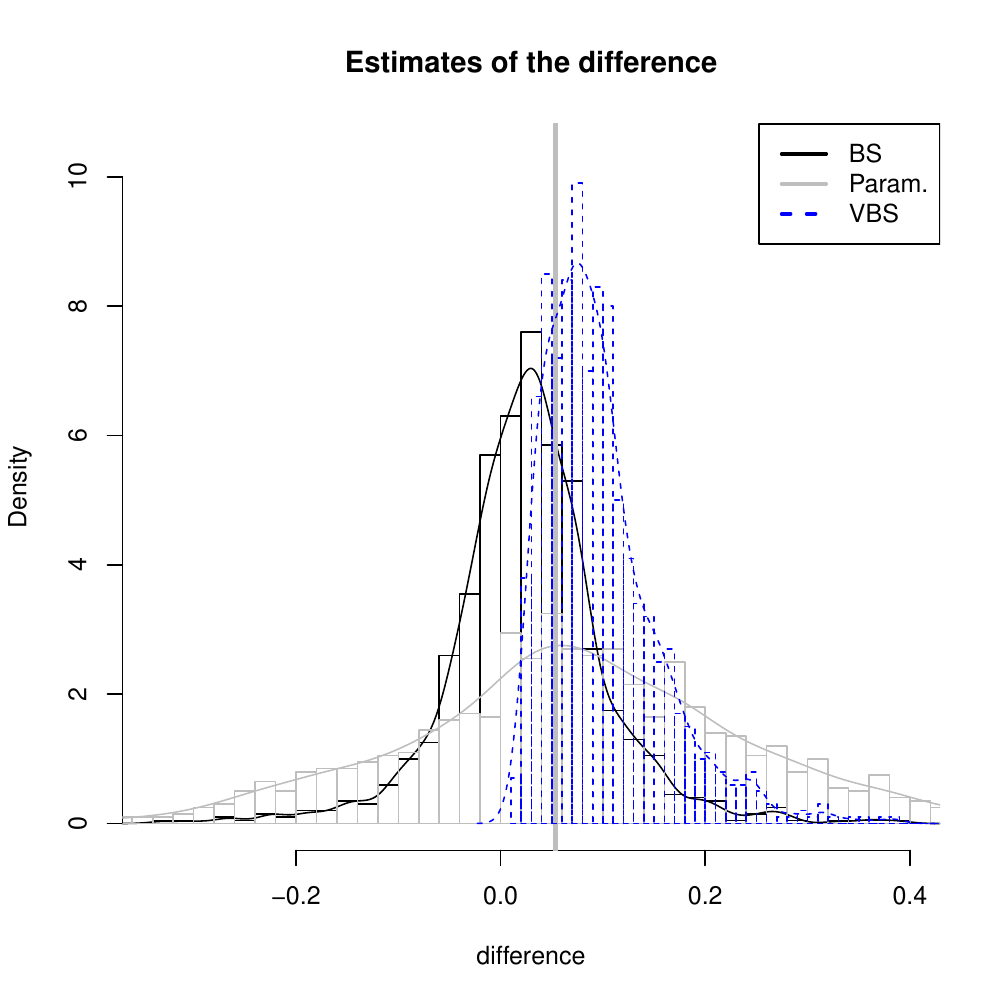}
	}
	\caption{\footnotesize Histograms of the estimates (Bootstrap, Parametric, Variance BS) of the asymptotic variance ${\bs \Sigma}_{PI}$ and the difference ${\bs \Sigma}_{LSE}-{\bs \Sigma}_{PI}$ (here the matrices are one dimensional).The true asymptotic value is illustrated by the vertical gray line. \label{fig:AV}}
\end{figure}

\clearpage

\subsection{Estimating the prediction error}\label{sec:pred_error}

The three estimates of ${\bold \Sigma}_{PI}$ can be used to estimate the prediction error in a standard fashion. Here, we describe how the variance bootstrap estimate together with the result of Proposition \ref{prop1}  provide an estimate of the excess risk.

By Proposition \ref{prop1}, the excess risk of $\hat{\boldsymbol \beta}_{LSE}$ can be approximated by
\[
E\Big[  Y -E(Y) - \sum_{j=1}^p {\beta}_j\{ { X}_j - E(X_j)\} \Big]^2
+ {\rm Trace} ({\bf M} {\boldsymbol \Sigma}_{LSE})=Var(Y)+ {\boldsymbol \beta}^T {\bf M} {\boldsymbol \beta} - 2\sum_{j=1}^p \beta_j Cov(Y,X_j) + {\rm Trace} ({\bf M} {\boldsymbol \Sigma}_{LSE}),
\]
and hence can be consistently estimated by
\[
\hat{\sigma}^2_Y + \hat{\boldsymbol \beta}_{LSE}^T \hat{\bf M} \hat{\boldsymbol \beta}_{LSE}-2\sum_{j=1}^p \{\hat{\boldsymbol \beta}_{LSE}\}_j  \overline{Cov}(Y,X_j) + {\rm Trace} \{\hat{\bf M} \hat{\boldsymbol \Sigma}_{LSE,BS}\},
\]
where $\hat{\sigma}^2_Y=\frac{1}{n} \sum_{i=1}^n (Y^{(i)} -\bar{Y})^2$, $\hat{\bf M}=\frac{1}{n+m} \sum_{i=1}^{n+m} ({\bf X}^{(i)}-\bar{\bf X})\left\{{\bf X}^{(i)}-\bar{\bf X}\right\}^{T}$ and $\overline{Cov}(Y,X_j)=\frac{1}{n}\sum_{i=1}^n(Y^{(i)} - \bar{Y})(X^{(i)}_{j}-\bar{X}_j)$ is the empirical covariance. Therefore, a consistent estimate to the excess risk ratio is 
\begin{equation}\label{eq:hat_ERR}
\widehat{ERR}_{PI}=\frac{\hat{\sigma}^2_Y + \hat{\boldsymbol \beta}_{LSE}^T \hat{\bf M} {\boldsymbol \beta}_{LSE}-2\sum_{j=1}^p \{\hat{\boldsymbol \beta}_{LSE}\}_j  \overline{Cov}(Y,X_j) + {\rm Trace} \{\hat{\bf M} (\hat{\boldsymbol \Sigma}_{LSE,BS}-\hat{\boldsymbol \Delta})\}}{\hat{\sigma}^2_Y + \hat{\boldsymbol \beta}_{LSE}^T \hat{\bf M} {\boldsymbol \beta}_{LSE}-2\sum_{j=1}^p \{\hat{\boldsymbol \beta}_{LSE}\}_j  \overline{Cov}(Y,X_j) + {\rm Trace} \{\hat{\bf M} \hat{\boldsymbol \Sigma}_{LSE,BS}\} },
\end{equation}
where $\hat{\boldsymbol \Delta}$ is defined by the Variance bootstrap.

\section{Sufficient and necessary conditions for improvement} \label{app:F}

{In this section we prove sufficient and necessary conditions under which $\sigma^2_{\rm diff} > 0$ in the one-dimensional case and ${\bs \Sigma}_{\rm diff} \ne {\bf 0}$ in the multidimensional case.

\subsection{One dimensional ${X}$}
Consider the simple linear model \eqref{eq:model1_1}. We have that $\hat{\beta}_{TI}$ and $\hat{\beta}_{PI}$ have smaller asymptotic variances if and only if $\sigma^2_{\rm diff}$ is not 0, where $\sigma^2_{\rm diff}$ is defined in Theorem \ref{thm:one_dim} (ii). The following proposition states a sufficient and necessary condition for $\sigma^2_{\rm diff} \ne 0$.

\begin{proposition}\label{prop:one}
	Under the conditions of Theorem \ref{thm:one_dim} (ii) we have that $\sigma^2_{\rm diff} \ne 0$ if and only if $\delta$ is uncorrelated with both $X^2$ and $X^3$
\end{proposition}

\noindent{\bf Proof:}
The total information estimator is the intercept estimator of the model \eqref{eq:modelXY}, which is 
\begin{equation*} 
W = \beta + a U_1  + b U_2  + \tilde{\delta}.  
\end{equation*}
Let ${\bf U}=(U_1,U_2)^T$ and ${\bf C}_U=E({\bf U}{\bf U}^T)$; notice that $E({\bf U})={\bf 0}$. We have that $\binom{a}{b} ={\bf C}_U^{-1}\binom{E(WU_1)}{E(WU_2)}$. Also,
\[
\sigma^2_{diff}=E\left({\delta X}-\tilde{\delta}\right)^2 = E\left\{(a-\alpha)U_1+(b-\beta)U_2 \right\}^2=\left( \begin{array}{c} a -\alpha \\ b -\beta \end{array}\right)^T {\bf C}_U\left( \begin{array}{c} a -\alpha \\ b -\beta \end{array}\right).
\]

Now let $\binom{\alpha'}{\beta'} = {\bf C}_U\binom{\alpha}{\beta}$, then $\binom{\alpha}{\beta} = {\bf C}_U^{-1}\binom{ \alpha'}{ \beta'}$ and therefore using the definition of $\binom{a}{b}$ it follows that 
\[
\sigma^2_{diff}=\left( \begin{array}{c} a -\alpha \\ b -\beta \end{array}\right)^T {\bf C}_U\left( \begin{array}{c} a -\alpha \\ b -\beta \end{array}\right)=
\left( \begin{array}{c} E(WU_1) -\alpha' \\ E(WU_2) -\beta' \end{array}\right)^T {\bf C}^{-1}_U\left( \begin{array}{c} E(WU_1) -\alpha' \\ E(WU_2) -\beta' \end{array}\right).
\]
Since, ${\bf C}_U$ is positive definite
$\sigma^2_{diff}=0$ if and only if $ E(WU_1) =\alpha'$ and $E(WU_2) = \beta'$.
Recall that $W=\frac{Y\{X-E(X)\}}{Var(X)}$, $U_1=\frac{X-E(X)}{Var(X)}$ and $U_2=\frac{\{X-E(X)\}X}{Var(X)} -1$. Therefore,
\[
E(WU_1)=E(U_1^2 Y)=E\{U_1^2(\alpha+\beta X + \delta) \}=\alpha E(U_1^2)+\beta E(U_1^2 X)+ E(U_1^2 \delta)=\alpha E(U_1^2)+\beta E(U_1 U_2)+ E(U_1^2 \delta)
\]
and $\alpha'=E(U_1^2)\alpha+E(U_1U_2) \beta$. It follows that $E(WU_1) -\alpha'=E(U_1^2 \delta)$. Similarly, $E(WU_2) -\beta'=E(U_1 U_2 \delta)$. Therefore, $\sigma^2_{diff}=0$ if and only if $E(U_1^2 \delta)=0$ and $E(U_1 U_2 \delta)=0$. The latter occurs if and only if $\delta$ is uncorrelated with both $X^2$ and $X^3$.
 \qed

\subsection{Multidimensional dimensional ${\bf X}$}
We now consider the general p-variate case \eqref{eq:model}. The following Proposition is parallel to the result for one dimensional $X$ and the proof is similar and hence omitted.

\begin{proposition}
	Under the conditions of Theorem \ref{thm:p_dim} (ii) we have that $\{{\bs \Sigma}_{\rm diff}\}_{j,j} \ne 0$ if and only if $\delta$ is uncorrelated with the vector $(X_{j\bullet}^2,X_{j\bullet}^2X_1,\ldots,X_{j\bullet}^2X_p)$.
\end{proposition}

\subsection{Binary covariates}\label{sec:binary}
In the one-dimensional model when $X$ is binary, then the linear model is trivially true. Therefore, $\delta$ is uncorrelated with both $X^2$ and $X^3$ and there is no improvement for $\hat{\beta}_{PI}$. In the multivariate case, when $X_j$ (say) is binary, then it seems intuitively that no improvement over the LSE is possible because there are no nonlinearities of $X_j$ in $\delta$. This intuitive argument is supported by many simulations (which are not reported here). However, we were not able to formally show that in the binary case $\delta$ is uncorrelated with the vector $(X_{j\bullet}^2,X_{j\bullet}^2X_1,\ldots,X_{j\bullet}^2X_p)$.    }

\end{document}